\documentclass[a4paper,11pt]{article}
\usepackage[english]{babel}
\usepackage{latexsym}
\usepackage{amssymb}
\usepackage{color}
\usepackage{graphicx}
\usepackage{multicol}
\usepackage{amsmath}
\usepackage[sans]{dsfont}
\usepackage{dsfont}
\usepackage{mathpazo} 
\usepackage[bb=boondox,scr=rsfso]{mathalfa}
\usepackage{mathrsfs}   
\usepackage{pifont}
\usepackage{booktabs}
\usepackage{enumerate}
\usepackage{hyperref}
\usepackage{url}
\usepackage{pxfonts}
\DeclareMathAlphabet{\mathpzc}{OT1}{pzc}{m}{it}

\definecolor{citecol}{rgb}{0,0,.7}
\definecolor{formulacol}{rgb}{.6,0,0}
\newenvironment{remark}{\par\noindent{\underline{\textsf{\textbf{Remark:}}}} }{\begin{flushright} $\blacktriangle$ \end{flushright}\par}

\newenvironment{proof}{\par\noindent{\textcolor{blue}{\underline{\textsc{\textbf{Proof}}}} }}{\textcolor{blue}{\begin{flushright} $\Box$ \end{flushright}}\par}
\newcommand{\EndTh}{\begin{flushright} \textcolor{blue}{$\mathbf{\square}$} \end{flushright}}

\newcommand{\bfun}{\left\{\begin{array}{ll}}
\newcommand{\efun}{\end{array}\right.}
\newcommand{\bfunB}{\begin{array}{ll}}
\newcommand{\efunB}{\end{array}}

\newcommand{\Prob}{\textsf{P}}
\newcommand{\kIdx}{k \in \text{K}}
\newcommand{\EmpProb}{\mathbb{P}}

\newcommand{\MyR}{\mathbb{R}}
\newcommand{\MyZ}{\mathbb{Z}}

\newcommand{\MyNN}{\mathbb{N}}

\newcommand{\nSeq}{\{1,\ldots,n\}}

\newcommand{\Exp}{\mathbb{E}}

\newcommand{\btheta}{\boldsymbol{\theta}}

\newcommand{\bz}{\mathbf{z}}
\newcommand{\bZ}{\mathbf{Z}}

\newcommand{\Beso}{\mathpzc{B}}
\newcommand{\Lp}{\textsf{L}}
\newcommand{\Ldue}{\textsf{L}^2}

\newcommand{\MyJ}{\mathrm{J}}
\newtheorem{teo}{Theorem}[section]
\newtheorem{cor}[teo]{Corollary}
\newtheorem{lem}[teo]{Lemma}

\renewcommand{\epsilon}{\varepsilon}

\renewcommand{\hat}{\widehat}                       

\newcommand{\bi}{\begin{itemize}}
\newcommand{\ei}{\end{itemize}}
\newcommand{\be}{\begin{eqnarray*}}
\newcommand{\ee}{\end{eqnarray*}}
\newcommand{\bequ}{\begin{equation}}
\newcommand{\eequ}{\end{equation}}
\title{A note on an Adaptive Goodness--of--Fit test with Finite Sample Validity for Random Design Regression Models}
\author{Pierpaolo Brutti\footnote{e--mail: \url{pierpaolo.brutti@uniroma1.it}} \\ \small{\textsf{Department of Statistics, Sapienza University of Rome}}}
\date{}
\begin{document}
\maketitle
\begin{center}
{\bf Abstract}
\begin{quote}
{\footnotesize
Given an i.i.d. sample $\{(X_i,Y_i)\}_{i\in\nSeq}$ from the random design regression model $Y = f(X) + \epsilon$ with $(X,Y) \in [0,1] \times [-M,M]$, in this paper we consider the problem of testing the (simple) null hypothesis ``$f = f_0$'', against the alternative ``$f \neq f_0$'' for a fixed $f_0 \in \Ldue([0,1],G_X)$, where $G_X(\cdot)$ denotes the marginal distribution of the design variable $X$. The procedure proposed is an adaptation to the regression setting of a multiple testing technique introduced by Fromont and Laurent \cite{Fromont:Laurent:2005}, and it amounts to consider a suitable collection of unbiased estimators of the $\Ldue$--distance $\textsf{d}_2(f,f_0) = \int {[f(x) - f_0 (x)]^2 \mathrm{d}G_X (x)}$, rejecting the null hypothesis when at least one of them is greater than its $(1-u_\alpha)$ quantile, with $u_\alpha$ calibrated to obtain a level--$\alpha$ test. To build these estimators, we will use the \emph{warped wavelet} basis introduced by Picard and Kerkyacharian \cite{KP04}.
We do not assume that the errors are normally distributed, and we do not assume that $X$ and $\epsilon$ are independent but, mainly for technical reasons, we will assume, as in most part of the current literature in learning theory, that $|f(x) - y|$ is uniformly bounded (almost everywhere). We show that our test is adaptive over a particular collection of approximation spaces linked to the classical Besov spaces.
\noindent
\medskip

\noindent{\em Keywords:}
Nonparametric Regression;
Random Design;
Goodness--of--fit;
Adaptive test;
Separation Rates;
Warped Wavelets;
U--statistics;
Multiple Test.
}
\end{quote}
\end{center}

\section{Introduction}\label{sec:CH2:Intro}
 
Consider the usual nonparametric regression problem with random design. In this model we observe an i.i.d. sample $\mathpzc{D}_n = \{\bZ_i = (X_i,Y_i)\}_{i\in\{1\ldots n\}}$ from the distribution of a vector $\bZ = (X,Y)$ where
\[
Y = f(X) + \epsilon,
\]
for $(X,\epsilon)$ a random vector with $\Exp(\epsilon | X) = 0$ and $\Exp(\epsilon^2|X) < \infty$ almost surely.
The regression function is known to belong to a subset $\mathcal{F}$ of $\Lp^2([0,1],G_X)$ for $G_X$ the marginal distribution of $X$. Let $f_0 \in \mathcal{F}$ be fixed. In this paper we consider the problem of testing the (simple) null hypothesis ``$\textsf{H}_0:f = f_0$'' against the alternative ``$\textsf{H}_1:f \neq f_0$''. Since $f \in \Ldue([0,1],G_X)$, it seems natural to consider a test statistic somehow linked to an estimator of the (weighted) $\Ldue$--distance $\textsf{d}_2(f,f_0) = \int {[f(x) - f_0 (x)]^2 \mathrm{d}G_X (x)}$. The approach considered in the present paper is an adaptation to the regression setting with random design of the work by Fromont and Laurent \cite{Fromont:Laurent:2005} for density models, and it amounts to consider a suitable collection of unbiased estimators for $\textsf{d}_2(f,f_0)$, rejecting the null hypothesis when at least one of them is greater than its $(1-u_\alpha)$ quantile, with $u_\alpha$ calibrated to obtain a level--$\alpha$ test.

After Ingster's seminal paper \cite{Ingster:1982}, and Hart's influential book \cite{Hart:1997}, many authors have been concerned with the construction of nonparametric tests on the unknown function that appears in a regression or Gaussian white noise model.
In papers like \cite{Ermakov:1991}, \cite{Ingster:1993a}, \cite{Hardle:Mammen:1993}, or more recently H\"{a}rdle and Kneip \cite{Hardle:Kneip:1999}, Lepski and Spokoiny \cite{Lepski:Spoko:1999}, Lepski and Tsybakov \cite{Lepski:Tsyba:2000}, and Gayraund and Pouet \cite{Gayraud:Pouet:2001}, the authors tackle the \emph{non--adaptive} case by specifying a particular functional/smoothness class to which $f(\cdot)$ belongs, and then evaluating the minimal separation/distance between the null hypothesis and the set of alternatives for which testing with a prescribed error probability is still possible.
Hard--coding the smoothness class of choice into any statistical procedure is clearly impractical and unattractive. For this reason, much of the effort has then be dedicated to explore the \emph{adaptive}, case where the smoothness level is also supposed to be unknown. So, for example, in\cite{Gayraud:Pouet:2005} and \cite{Horo:Spoko:2001} Gayraund and Pouet on one side and Horowitz and Spokoiny on the other deal with the adaptive case for a composite null hypothesis and suitable smoothness classes (e.g. H\"{o}lder spaces), whereas Fromont and L\'{e}vy--Leduc  in \cite{Fromont:Leduc:2006} consider the problem of periodic signal detection in a Gaussian fixed design regression framework, when the signal belongs to some periodic Sobolev balls. 
Fan, Zhang and Zhang \cite{Fan:Zhang:Zhang:2001} and Fan and Zhang \cite{Fan:Zhang:2004}, using a generalized likelihood ratio, give adaptive results when the alternatives lie in a range of Sobolev ball, also highlighting a nonparametric counterpart of the so called Wilks phenomenon well known in parametric inference.
In \cite{Spokoiny:1996} and \cite{Spokoiny:1998} Spokoiny considers testing a simple hypothesis under a Gaussian white noise model over an appropriate collection of Besov balls.
Quite relevant is also the work of Baraud, Huet and Laurent \cite{Baraud:Huet:Laurent:2003a} where the assumption on $f(\cdot)$ are reduced to a minimum thanks to the adoption of a discrete distance that approximate the usual $\Lp^2$--norm to measure separation between the null and the alternative hypothesis.

Similar problems have been widely studied in the testing literature. To briefly summarize the basic notions and notation regarding hypothesis testing, consider the following general setting where we have an observation $Y$ coming from a distribution $G_Y(\cdot)$, and we are interested in testing the (composite) null hypothesis $\textsf{H}_0: G_Y \in \mathscr{G}_0$, where $\mathscr{G}_0$ denotes a families of probability measures, against the alternative $\textsf{H}_1:G_Y  \notin \mathscr{G}_0$. To accomplish this task, we need to define a test function $\textsf{T}_\alpha(Y)$; that is, a measurable function of $Y$ that takes values in $\{0,1\}$, such that, given a testing level $\alpha \in (0,1)$, we reject $\textsf{H}_0$ every time $\textsf{T}_\alpha(Y) = 1$. The value $\alpha$ is the testing level of our procedure in the sense that we require
\[
\mathop {\sup }\limits_{G  \in \mathscr{G}_0 } \Prob_{G } \big\{ {\textsf{T}_\alpha  (Y) = 1} \big\} \leqslant \alpha.
\]
For each $G \notin \mathscr{G}_0$, the Type II error of our testing procedure on $G(\cdot)$, is defined by
\[
\beta \big( {G,\textsf{T}_\alpha  (Y)} \big) = \Prob_G \big\{ {\textsf{T}_\alpha  (Y) = 0} \big\},
\]
whereas the power of the test on $G(\cdot)$ is given by
\[
\pi \big({G,\textsf{T}_\alpha(Y)} \big) = 1 - \beta \big({G,\textsf{T}_\alpha  (Y)} \big).
\]
Of course, an easy way to choose a testing method would be to select the most powerful one (i.e. the one with smaller Type II error) within the class of level--$\alpha$ tests. In general, the closer $G \notin \mathscr{G}_0$ is to $\mathscr{G}_0$, the more difficult is to \emph{separate} the null from the alternative hypothesis, and consequently, the smaller is the power of a test on that particular $G(\cdot)$. This obvious fact naturally leads to define the notion of \emph{separation rate} of a test $\textsf{T}_\alpha(Y)$ over a functional class $\mathscr{G}$, with respect to a distance $\textsf{d}(\cdot)$, as follow
\[
\rho \left( {\textsf{T}_\alpha (Y),\mathscr{G},\beta } \right) = \inf \left\{ {\rho  > 0:\mathop {\sup }\limits_{G \in \mathscr{G}:\textsf{d}(G,\mathscr{G}_0 ) \geqslant \rho } \Prob_G \big\{ {\textsf{T}_\alpha (Y) = 0} \big\} \leqslant \beta } \right\}.
\]
In words, $\rho \left( {\textsf{T}_\alpha (Y),\mathscr{G},\beta } \right)$ is the minimal distance from $\mathscr{G}_0$ starting from which our testing procedure has a Type II error smaller than $\beta$ uniformly over $\mathscr{G}$. From here, we can immediately define the (non--asymptotic) $(\alpha,\beta)$ \emph{minimax rate} of testing over the class $\mathscr{G}$ as follow
\[
\mathop {\inf }\limits_{\textsf{T} \in \mathcal{T}_\alpha  } \rho \left( {\textsf{T},\mathscr{G},\beta } \right),
\]
where $\mathcal{T}_\alpha$ denotes the class of test statistics that are associated to $\alpha$--level tests.

If we have a complete characterization of the class $\mathscr{G}$, in general we are able to build a testing procedure that explicitly depends on $\mathscr{G}$, and attains the minimax separation rate over $\mathscr{G}$ itself. However, as we already said, it is extremely unsatisfying to have $\mathscr{G}$ hard--coded in our technique. A more interesting task, in fact, would be to build adaptive testing methods that simultaneously (nearly) attain the minimax separation rate over a broad range of reasonable classes without using any prior knowledge about the law of the observations. Eubank and Hart \cite{Eubank:Hart:1992} propose to test that a regression function is identically zero using a test function based on the Mallow's $C_p$ penalty. Antoniadis, Gijbels and Gr\'{e}goire \cite{Anto:Gije:Grego:1997}, once again in a regression setting, develop an automatic model selection procedure that they also apply to the same testing problem. Spokoiny \cite{Spokoiny:1996}, instead, considers a Gaussian white noise model $\mathrm{d}X(t) = f(t)\mathrm{d}t + \epsilon \mathrm{d}W(t)$, and propose to test ``$f \equiv 0$'' adaptively using a wavelet based procedure. He also study the (asymptotic) properties of his approach and show that, in general, adaptation is not possible without some loss of efficiency of the order of an extra $\log\log(n)$ factor, where $n$ is the sample size (see Section \ref{subsec:CH2:TestUniSepRates}). In the same setting, Ingster \cite{Ingster:2000} builds an adaptive test based on chi--square statistics, and study its asymptotic properties.

Many authors have also considered the problem of testing convex or \emph{qualitative} hypothesis like the monotonicity of the regression function: Bowman, Jones and Gijbels \cite{Bow:Jones:Gij:1998}; Hall and Heckman \cite{Hall:Heckman:2000}; Gijbels, Hall, Jones and Koch \cite{Gij:Hall:Jones:Koch:2000}; Ghosal, Sen and van der Vaart \cite{Ghosal:Sen:vdVaart:2001}, are just a few examples. In \cite{Dumb:Spoko:2001}, instead, D\"{u}mbgen and Spokoiny consider the problem of testing the positivity, monotonicity and convexity of the function $f(\cdot)$ that appears in a Gaussian white noise model. They also evaluate the separation rates of their procedures showing in this way their optimality. See also Juditsky and Nemirovski \cite{Judi:Nemirovski:2002}.

The literature regarding goodness--of--fit testing in a density model is also vast. Bickel and Ritov \cite{Bickel:Ritov:1992}; Ledwina \cite{Ledwina:1994}; Kallenberg and Ledwina \cite{Kallenberg:Ledwina:1995}; Inglot and Ledwina \cite{Inglot:Ledwina:1996}; Kallenberg \cite{Kallenberg:2002}, for instance, propose tests inspired by Neyman \cite{Neyman:1937} where the parameter that enter the definition of the test statistic (in general a smoothing parameter) is estimated by some automatic data dependent criterion like \textsf{BIC}. In general, only the asymptotic performances of these tests have been studied in some detail. 

The pre--testing approach considered in the paper by Fromont and Laurent \cite{Fromont:Laurent:2005} has been initiated by Baraud, Huet and Laurent \cite{Baraud:Huet:Laurent:2001,Baraud:Huet:Laurent:2003a,Baraud:Huet:Laurent:2003b} for the problem of testing linear or qualitative hypotheses in the Gaussian regression model. One nice feature of their approach is that the properties of the procedures are non asymptotic. For any given sample size $n$, the tests have the desired level and we are able to characterize a set of alternatives over which they have a prescribed power. It is interesting to notice that the method proposed by Fromont and Laurent to build a test function essentially amounts to penalize by the appropriate quantile under the null, an unbiased estimator of projections of the $\Ldue$--distance between densities.
Other papers where U--statistics have been used to build test functions are: \cite{Dette:Neumeyer:2001,Neumeyer:Dette:2003,Munk:Dette:1998,Butucea:Tribouley:2006}.

This paper is organized as follow. In Section \ref{sec:CH2:GoodOfFit} we describe the testing procedure. In Section \ref{subsec:CH2:TestPower} we review the concept of \emph{warped wavelet basis} proposed in \cite{KP04} and we establish the type of alternatives against which our test has a guaranteed power. Then, in Section \ref{subsec:CH2:TestUniSepRates} together with a brief simulation study, we show that our procedure is adaptive over some collection of \emph{warped} Besov spaces in the sense that it achieves the optimal ``adaptive'' rate of testing over all the members of this collection simultaneously. Finally Section \ref{sec:CH2:ProofsTest} contains the proofs of the results presented.

\section{A Goodness--of--Fit Test}\label{sec:CH2:GoodOfFit}

As anticipated in the previous section, the framework we shall work with in this paper is the usual nonparametric regression problem with random design. In this model we observe an i.i.d. sample $\mathpzc{D}_n = \{\bZ_i = (X_i,Y_i)\}_{i\in\{1\ldots n\}}$ from the distribution of a vector $\bZ = (X,Y)$ described structurally as
\[
    Y = f(X) + \epsilon,
\]
for $(X,\epsilon)$ a random vector with $\Exp(\epsilon | X) = 0$ and $\Exp(\epsilon^2|X) < \infty$ almost surely.
The regression function is known to belong to a subset $\mathcal{F}$ of $\Lp^2([0,1],G_X)$ for $G_X$ the marginal distribution of $X$, which is assumed known. As explained in \cite{Robins:vdVaart:2004}, the assumption on $G_X$ is surely unpleasant but unavoidable: the radius of the confidence set will be inflated in varying amount depending on the conditions imposed on $G_X$ so we postpone the treatment of this case to a forthcoming paper.
The variance function $\sigma^2(x) = \Exp(\epsilon^2|X=x)$ need not to be known, although a known upper bound on $\|\sigma^2\|_{\infty}$ is needed. We do not assume that the errors are normally distributed, and we do not assume that $X$ and $\epsilon$ are independent but, mainly for technical reasons, we will assume, as in most part of the current literature in learning theory (see \cite{CS}), that $|f(x) - y|$ is uniformly bounded (almost everywhere) by a positive constant $M$. Doing so, all the proofs will be greatly simplified without moving too far away from a realistic (although surely not minimal) set of assumptions (in particular considering the finite--sample scope of the analysis). Clearly this condition overrules the one on the conditional variance mentioned before.
        
As it is often the case in nonparametric statistics, we could cast this example into a problem of estimating a sequence $\btheta = [\theta_1,\theta_2,\ldots] \in \ell^2$ of parameters by expanding $f(\cdot)$ on a fixed orthonormal basis $\{e_\ell\}_{\ell\in\MyNN}$ of $\Lp^2([0,1],G_X)$. The Fourier coefficients take the form 
\[
    \theta_\ell = \langle f, e_\ell \rangle_{\Lp^2(G_X)} = \Exp_{(X,Y)}[Y \, e_\ell(X)],
\]
and they can be estimated unbiasedly by $W_\ell = \tfrac{1}{n}\sum\nolimits_{i = 1}^n {Y_i \;e_\ell (X_i )}$, although it appears not so useful to move directly in sequence space by considering $[W_1,W_2,\ldots]$ as the observation vector. What we propose is a goodness--of--fit test similar to the one introduced in \cite{Fromont:Laurent:2005}. To describe it, let $f_0(\cdot)$ be some fixed function in $\Lp^2([0,1],G_X)$ and $\alpha \in (0,1)$. Now suppose that our goal is to build a level--$\alpha$ test of the null hypothesis $\textsf{H}_0: f \equiv f_0$ against the alternative $\textsf{H}_1: f \neq f_0$ from the data $\{\bZ_i\}_{i\in\nSeq}$. The test is based on an estimation of 
\[
\|f - f_0\|^2_{\Lp^2(G_X)} = \|f\|^2_{\Lp^2(G_X)} + \|f_0\|^2_{\Lp^2(G_X)} - 2 \langle f, f_0 \rangle_{\Lp^2(G_X)}.
\]
Since the last (linear) term $\langle f, f_0 \rangle_{\Lp^2(G_X)}$ can be easily estimated by the empirical estimator $\tfrac{1}{n}\sum\nolimits_{i = 1}^n {Y_i \;f_0 (X_i )}$, the key problem is the estimation of the first term $\|f\|^2_{\Lp^2(G_X)}$.
Adapting the arguments in \cite{Laurent:2005a}, we can consider an at most countable collection of linear subspaces of $\Lp^2([0,1],G_X)$ denoted by $\mathcal{S} = \{\text{S}_k\}_{k \in \text{K}}$. For all $k\in\text{K}$, let $\{e_\ell\}_{\ell \in \mathcal{I}_k}$ be some orthonormal basis of $\text{S}_k$. The estimator
\begin{equation}\label{eq:CH2:UStat}
\hat \theta _{n,k}  = \frac{1}{{n(n - 1)}}\sum\limits_{i = 2}^n {\sum\limits_{j = 1}^{n - 1} {\left[ {\sum\limits_{\ell  \in \mathcal{I}_k } {\big\{ {Y_i e_\ell  (X_i )} \big\} \cdot \big\{ {Y_j e_\ell (X_j )} \big\}} } \right]} }  = \frac{1}{{n(n - 1)}}\sum\limits_{i = 2}^n {\sum\limits_{j = 1}^{n - 1} {h_k(\bZ_i ,\bZ_j )} } ,
\end{equation}
is a U--statistic of order two for $\|\Pi_{\text{S}_k}(f)\|^2_{\Lp^2(G_X)}$ -- where $\Pi_{\text{S}_k}(\cdot)$ denotes the orthogonal projection onto $\text{S}_k$ -- with kernel
\[
h_k(\bz_1,\bz_2) = {\sum\limits_{\ell  \in \mathcal{I}_k } {\left\{ {y_1 e_\ell  (x_1 )} \right\} \cdot \left\{ {y_2 e_\ell (x_2 )} \right\}} },\quad \bz_i = (x_i,y_i),\;i\in\{1,2\}.
\]
Then, for any $\kIdx$, $ \|f - f_0\|^2_{\Lp^2(G_X)}$ can be estimated by
\begin{eqnarray}\label{eq:CH2:riskDecompo}
\widehat{\textsf{R}}_{n,k}  &=& \widehat\theta_{n,k} + \| {f_0 } \|_{\Lp^2(G_X )}^2  - \frac{2}
                           {n}\sum_{i = 1}^n {Y_i \;f_0 (X_i )} = \hfill \nonumber \\
&{\mathop  = \limits^{\left( \diamondsuit  \right)}}& 
    \widetilde{U}_{n,k} + 2(\EmpProb_n  - \Prob) \big( {\Pi _{\text{S}_k } (f) - f} \big) - 
                            \| {\Pi _{\text{S}_k } (f) - f} \|_{\Lp^2 (G_X )}^2  +  \hfill \\
                        & & + 2(\EmpProb_n  - \Prob)\big( {f - f_0} \big) + \| {f - f_0 } \|_{\Lp^2 (G_X )}^2, \nonumber
\end{eqnarray}
where ,
\bequ\label{eq:CH2:2ndOrderUstat}
\widetilde{U}_{n,k}  = \frac{1}{{n(n - 1)}}\sum\limits_{i = 2}^n {\sum\limits_{j = 1}^{n - 1} {\left[ {\sum\limits_{\ell  \in \mathcal{I}_k } {\big\{ {Y_i e_\ell  (X_i ) - \theta _\ell  } \big\} \cdot \big\{ {Y_j e_\ell (X_j ) - \theta _\ell  } \big\}} } \right]} }  = \frac{1}
{{n(n - 1)}}\sum\limits_{i = 2}^n {\sum\limits_{j = 1}^{n - 1} {g_{k} (\bZ_i ,\bZ_j )} } ,
\eequ
and, for each $w\in \Ldue$,
\[
\EmpProb_n(w) = \frac{1}{n}\sum\limits_{i = 1}^n {Y_i w(X_i )} = \left\langle {f,w}\right\rangle _n
\quad \text{and} \quad \
\Prob(w) = \int {f(x)w(x)\mathrm{d}G_X (x)}  = \left\langle {f,w} \right\rangle _{\Ldue(G_X )},
\]
so that
\[
\Exp_{(X,Y)} \left\{ {\EmpProb_n(w)} \right\} = \frac{1}{n}\,n\,\Exp_{(X,Y)} \left\{ {Y \,w(X)} \right\} = \Exp_{G_X } \left\{ {f(X)\,w(X)} \right\} = \int {f(x)\,w(x)\mathrm{d}G_X (x)}  = \Prob(w).
\]
The equality $\left( \diamondsuit  \right)$ can be derived from the \emph{Hoeffding decomposition} of $\hat{\theta}_{n,k}$ as explained in Section \ref{sec:CH2:ProofsTest}.

Now that we have an estimator $\widehat{\textsf{R}}_{n,k}$, lets denote by $r_{n,k}(u)$ its $1-u$ quantile under $\textsf{H}_0$, and consider
\[
u_\alpha   = \sup \left\{ {u \in (0,1):\Prob^{\otimes n}_{f_0 } \left[ {\mathop {\sup }\limits_{k \in K} \big\{ {\widehat{\textsf{R}}_{n,k} - r_{n,k} (u)} \big\} > 0} \right] \leqslant \alpha } \right\},
\]
where $\Prob^{\otimes n}_{f_0 }\{\cdot\}$ is the law of the observations $\{\bZ_i\}_{i\in\nSeq}$ under the the null hypothesis. Then introduce the test statistics $\textsf{R}_\alpha$ defined by
\[
\textsf{R}_\alpha   = \mathop {\sup }\limits_{k \in K} \big\{ {\widehat{\textsf{R}}_{n,k}  - r_{n,k} (u_\alpha  )} \big\},
\]
so that we reject the null whenever $\textsf{R}_\alpha$ is positive.

This method, adapted to the regression setting from \cite{Fromont:Laurent:2005}, amounts to a multiple testing procedure. Indeed, for all $k \in \text{K}$, we construct a level--$u_\alpha$ test by rejecting $\textsf{H}_0: f \equiv f_0$ if $\widehat{\textsf{R}}_{n,k}$ is greater than its $(1-u_\alpha)$ quantile under $\textsf{H}_0$. After this, we are left with a collection of tests and we decide to reject $\textsf{H}_0$ if, for some of the tests in the collection, the hypothesis is rejected. In practice, the value of $u_\alpha$ and the quantile $\{r_{n,k}(u_\alpha)\}_{\kIdx}$ are to be estimated by (smoothed) bootstrap (see \cite{Hormann:Ley:Derf:2004,Shao:Tu:1995}).

\subsection{Power of the Test}\label{subsec:CH2:TestPower}

Both the practical and theoretical performances of the proposed test, depend strongly on the orthogonal system we adopt to generate the collection of linear subspaces $\{S_k\}_{\kIdx}$. In dealing with a density model, Fromont and Laurent \cite{Fromont:Laurent:2005}, consider a collection obtained by mixing spaces generated by constants piecewise functions (Haar basis), scaling functions from a wavelet basis, and, in the case of compactly supported densities, trigonometric polynomial. Clearly these bases are not orthonormal in our weighed space $\Ldue([0,1],G_X)$, hence we have to consider other options. 

The first possibility that comes to mind is to use one of the usual wavelet bases since, as proved by Haroske and Triebel in \cite{Haroske:Triebel:2005} (see also \cite{Garcia:Martell:2001} and \cite{Kuhn:etal:2005}), these systems continue to be unconditional \emph{Schauder} bases for a whole family of weighted Besov spaces once we put some polynomial restriction on the growth of the weight function. 

Although appealing, this approach has some evident drawbacks once applied to our setting from a theoretical (we must impose some counterintuitive conditions on the marginal $G_X(\cdot)$), as well as practical (we can not use the well--known \emph{fast wavelet transform} anymore, see \cite{Brutti:2005b}) point of view.

As one can see looking at the proofs of Section \ref{sec:CH2:ProofsTest}, a basis that proved to fit perfectly in the present framework, is the so--called \emph{warped wavelet} basis studied by Kerkyacharian and Picard in \cite{KP04,Kerk:Picard:2005}. The idea is as follow. For a signal observed at some design points, $Y(t_i)$, $i \in \{1,\ldots,2^\mathrm{J}\}$, if the design is regular ($t_k = k/2^\mathrm{J}$), the standard wavelet decomposition algorithm starts with $s_{\mathrm{J},k} = 2^{\mathrm{J}/2}Y(k/2^{\mathrm{J}})$ which approximates the scaling coefficient $\int Y(x)\phi_{\mathrm{J},k}(x)\mathrm{d}x$, with $\phi_{\mathrm{J},k}(x) = 2^{\mathrm{J}/2}\phi(2^\mathrm{J} x - k)$ and $\phi(\cdot)$ the so--called scaling function or father wavelet (see \cite{Mallat:1998} for further information). Then the cascade algorithm is employed to obtain the wavelet coefficients $d_{j,k}$ for $j \leqslant \mathrm{J}$, which in turn are thresholded. If the design is not regular, and we still employ the \emph{same} algorithm, then for a function $H(\cdot)$ such that $H(k/2^\mathrm{J}) = t_k$, we have $s_{\mathrm{J},k} = 2^{\mathrm{J}/2}Y(H(k/2^{\mathrm{J}}))$. Essentially what we are doing is to decompose, with respect to a standard wavelet basis, the function $Y(H(x))$ or, if $G \circ H(x) \equiv x$, the original function $Y(x)$ itself but with respect to a new \emph{warped} basis $\{\psi_{j,k}(G(\cdot))\}_{(j,k)}$.

In the regression setting, this means replacing the standard wavelet expansion of the function $f(\cdot)$ by its expansion on the new basis $\{\psi_{j,k}(G(\cdot))\}_{(j,k)}$, where $G(\cdot)$ is adapting to the design: it may be the distribution function of the design $G_X(\cdot)$, or its estimation when it is unknown (not our case). An appealing feature of this method is that it does not need a new algorithm to be implemented: just standard and widespread tools (we will use this nice feature of the warped bases in the companion paper \cite{Brutti:2005b}).

It is important to notice that a warped wavelet basis is, automatically, an orthonormal system in $\Ldue([0,1],G_X)$. In fact, if, for easy of notation, we index the basis functions by mean of the set $\mathscr{D} \equiv \mathscr{D}([0,1])$ of dyadic cubes of $\MyR$ contained in $[0,1]$, i.e. we set $\psi_{j,k}(\cdot) \equiv \psi_{\textsf{I}}(\cdot)$, then for each $\textsf{I}_1$, $\textsf{I}_2$ in $\mathscr{D}$, we have
\[
\left\langle {\psi _{\textsf{I}_1 } \circ G_X,\psi _{\textsf{I}_2 } \circ G_X} \right\rangle _{L^2 (G_X )}  = \int {\psi _{\textsf{I}_1 } (G_X (x))\psi _{\textsf{I}_2 } (G_X (x))\mathrm{d}G_X (x)} = \int {\psi _{\textsf{I}_1 } (y)\psi _{\textsf{I}_2 } (y)\mathrm{d}y}  = \delta _{\textsf{I}_1 ,\textsf{I}_2 },
\]
where the last equality comes from the fact that we can build our warped basis from a (boundary corrected) wavelet system, orthonormal with respect to the Lebesgue measure in $[0,1]$ (see \cite{Cohen:Daub:Vial:1993,Johnstone:Silver:2004}, and Chapter $7.5$ in \cite{Mallat:1998}).

Now, to extract a basis out of a warped system, we surely need to impose restrictions on the design distribution $G_X(\cdot)$. As a matter of fact, it is easy to prove that the orthonormal system $\{\psi_\textsf{I}(G)\}_{\textsf{I}\in\mathscr{D}}$ or, equivalently, the system $\{ \phi _{\mathrm{J},k} (G)\} _{k \in \{1,\ldots,\overline{k}(\mathrm{J})\}}$ of scaling functions at any fixed resolution level $\mathrm{J}$, is \emph{total} in $\Ldue([0,1],G_X)$ if $G_X(\cdot)$ is absolutely continuous with density $g_X(\cdot)$ -- with respect to the Lebesgue measure -- bounded from below and above. Of course, this condition is only sufficient and unnecessarily stringent but also simple enough to fit perfectly our desiderata.

At this point, for each $\mathrm{J}$, we have a system of scaling functions $\{\phi_{\mathrm{J},k}(G)\}_k$ orthonormal in $\Ldue([0,1],G_X)$ that we can use to generate the subspaces $\mathcal{S} = \{\text{S}_{\mathrm{J}}\}_{\mathrm{J}\in\MyNN}$ where we have slightly changed the indexing notation: from $k$ to $\mathrm{J}$. So let
\[
\text{S}_\mathrm{J} = \textsf{span}\big\{ \{\phi_{\mathrm{J},k}(G)\}_{k\in\MyZ} \big\} \quad \text{with} \quad \mathrm{J} \in \big\{0,\ldots,\overline{\mathrm{J}}(n)\big\} \triangleq \mathcal{J}_n,
\]
and
\[
\hat{\theta}_{n,\MyJ}  = \frac{1}{{n(n - 1)}}\sum\limits_{i = 2}^n {\sum\limits_{j = 1}^{n - 1} {\left[ {\sum\limits_{k = 1}^{\overline k (\MyJ)} {\left\{ {Y_i \phi _{\MyJ,k} (G(X_i ))} \right\} \cdot \left\{ {Y_j \phi _{\MyJ,k} (G(X_j ))} \right\}} } \right]} }  = \frac{1}
{{n(n - 1)}}\sum\limits_{i = 2}^n {\sum\limits_{j = 1}^{n - 1} {h_{\MyJ} (\bZ_i, \bZ_j )} }.
\]
For all $\MyJ \in \mathcal{J}_n$, we set
\[
\widehat{\textsf{R}}_{n,\MyJ}  = \widehat\theta_{n,\MyJ} + \| {f_0 } \|_{\Lp^2(G_X )}^2  - \frac{2}{n}\sum_{i = 1}^n {Y_i \;f_0 (X_i )}.
\]
The test statistic we consider is
\bequ\label{eq:CH2:TestStat}
\textsf{R}_\alpha   = \mathop {\sup }\limits_{\MyJ \in \mathcal{J}_n} \big\{ {\widehat{\textsf{R}}_{n,\MyJ} - r_{n,\MyJ} (u_\alpha  )} \big\},
\eequ
where $r_{n,\MyJ} (u_\alpha)$ is defined in Section \ref{sec:CH2:GoodOfFit}.

The following theorem, which mimics Theorem $1$ in \cite{Fromont:Laurent:2005}, describes the class of alternatives over which the test has a prescribed power. The proof can be found in Section \ref{sec:CH2:ProofsTest}.
\begin{teo}\label{teo:CH2:PowerTest}
Let $\{\bZ_i = (X_i,Y_i)\}_{i\in\nSeq}$ be an i.i.d. sequence from the distribution of a vector $\bZ = (X,Y)$ described structurally by the nonparametric regression model
\[
Y = f(X) + \epsilon,
\]
for $(X,\epsilon)$ a random vector with $\Exp(\epsilon|X)=0$ and $\Exp(\epsilon^2|X)<+\infty$. Assume further that $f_0(\cdot)$ and the unknown regression function $f(\cdot)$ belong to $\Ldue([0,1],G_X)$ for $G_X(\cdot)$ the marginal distribution of $X$, assumed know and absolutely continuous with density $g_X(\cdot)$ bounded from below and above. Finally assume that $|f(x)-y|$ is uniformly bounded (almost everywhere) by a positive constant $M$.

Now let $\beta \in (0,1)$. For all $\gamma \in (0,2)$, there exist positive constants $C_1 \equiv C_1(\beta)$ and $C_2 \equiv C_2(\beta,\gamma,\tau_\infty,M,\|f_0\|_\infty)$ such that, defining
\[
\textsf{V}_{n,\MyJ} (\beta ) = \frac{C_1}{n}\left\{ {\tau_\infty \cdot \sqrt {2^\MyJ } + \frac{M^2}{n}2^\MyJ } \right\} + \frac{C_2}{n},
\]
with $\tau_\infty = \|f\|^2_\infty + \|\sigma^2\|_\infty$, then, for every $f(\cdot)$ such that
\[
\| {f - f_0 } \|_{\Ldue(G_X )}^2  > (1 + \gamma )\mathop {\inf }\limits_{\MyJ \in \mathcal{J}_n } \left\{ {\| {f - \Pi_{\text{S}_\MyJ } (f)} \|_{\Ldue(G_X )}^2  + r_{n,\MyJ} (u_\alpha  ) + \textsf{V}_{n,\MyJ} (\beta )} \right\},
\]
the following inequality holds:
\[
\Prob^{\, \otimes n}_f\big\{\textsf{R}_\alpha \leqslant 0\big\} \leqslant \beta.
\]
\end{teo}

\subsection{Uniform Separation Rates}\label{subsec:CH2:TestUniSepRates}

Now that we know against what kind of alternatives our multiple testing procedure has guaranteed power, we can move on, and examine the problem of establishing uniform separation rates over well--suited functional classes included in $\Ldue([0,1],G_X)$. We will start by defining for all $s>0$, $R>0$, and $M>0$, the following (linear) approximation space (see the review by DeVore \cite{DeVore:1998}):
\bequ\label{eq:CH2:ApprSpace}
\mathpzc{A}^s(R,M,G_X) = \big\{w\in\Ldue([0,1],G_X):\|w\|_\infty \leqslant M, \text{ and } \|w - \Pi_{\text{S}_\MyJ}(w)\|^2_{\Ldue(G_X)} \leqslant R^2 2^{-2 \MyJ \, s} \big\}.
\eequ
When $\mathrm{d}G_X(x) = \mathrm{d}x$ is the Lebesgue measure, $\mathpzc{A}^s(R,M,\mathrm{d}x)$ is strictly related to the following Besov body
\[
\Beso^{2,s}_\infty(R) = \left\{ w \in \Ldue(\mathrm{d}x,\MyR): \sum\limits_{k \in \MyZ} {d_{\MyJ,k}^2 } \leqslant R^2 2^{ - 2 \MyJ \, s} \right\}, \quad \text{with} \quad d_{\MyJ,k} = \left\langle {w,\psi _{\MyJ,k} } \right\rangle _{\Ldue (\MyR)},
\]
since
\[
\Beso_\infty ^{2,s} (R) \cap \left\{ {w: \| w \|_\infty   \leqslant  M } \right\} \subset \mathpzc{A}^s \left( {\tfrac{R}{{\sqrt {1 - 4^{ - s} } }},M,\mathrm{d}x} \right).
\]
In our case, instead, it is a bit less clear how to ``visualize'' the content of $\mathpzc{A}^s(R,M,G_X)$ in terms of common smoothness classes like Besov, H\"{o}lder or Sobolev body that admit alternative definitions in terms of geometric quantities like the \emph{modulus of smoothness} (see \cite{Devore:Lorentz:1993}).
The easiest way, probably, is to notice that, for each $w \in \Ldue([0,1],G_X)$
\[
\left\| {w - \Pi _{S_\MyJ } (w)} \right\|_{\Ldue(G_X )}^2  = \left\| {w(G_X^{ - 1} ) - \Pi _{S_\MyJ } \big(w(G_X^{ - 1} )\big)} \right\|_{\Ldue(\mathrm{d}x)}^2, 
\]
where the norm in the right hand side is taken with respect to the Lebesgue measure and
\[
G_X^{-1}(x) = \inf\{t \in \MyR : G_X(t) \geqslant x\}
\]
is the \emph{quantile function} of the design distribution $G_X(\cdot)$. Consequently,
\[
f \in \mathpzc{A}^s(R,M,G_X) \Leftrightarrow f(G_X^{-1}) \in \mathpzc{A}^s(R,M,\mathrm{d}x) \supset
\Beso^{2,s}_\infty \big(R\sqrt{1-4^{-s}}\big) \cap \{f:\|f\|_\infty \leqslant M\},
\]
so that the regularity conditions that hide behind the definition of the approximation space $\mathpzc{A}^s(R,M,G_X)$ could be expressed more explicitly in terms of the \emph{warped} function $f\circ G^{-1}_X(\cdot)$, mixing the smoothness of $f(\cdot)$ with the (very regular, indeed) design $G_X(\cdot)$ (see \cite{KP04,Kerk:Picard:2005} for further information and discussions).

The following corollary gives upper bounds for the uniform separation rates of our procedure over the class $\mathpzc{A}^s(R,M,G_X)$.

\begin{cor}\label{cor:CH2:unifSepRates}
Let $\textsf{R}_\alpha$ be the test statistic defined in Equation (\ref{eq:CH2:TestStat}). Assume that $n \geqslant 16$, and $\MyJ \in \mathcal{J}_n = \{0,\ldots,\overline{\mathrm{J}}(n)\}$ with
\[
2^{\overline{\mathrm{J}}(n)}  = \tfrac{n^2}{[\log\log(n)]^3}.
\]
Let $\beta \in (0,1)$. For all $s>0$, $M>0$, and $R>0$, there exist some positive constant $C \equiv C(s,\alpha,\beta,M,\|f_0\|_\infty)$ such that if $f \in \mathpzc{A}^s(R,M,G_X)$ and satisfies
\[
\| {f - f_0 } \|_{\Ldue(G_X )}^2  > C\left\{ {R^{\tfrac{2}{{4s + 1}}} 
    \left[ {\frac{{\sqrt {\log \log (n)} }}{n}} \right]^{\tfrac{4s}{{4s + 1}}}  + 
    R^2 \left[ {\frac{{(\log \log (n))^3 }}{{n^2 }}} \right]^{2s}  + 
    \frac{{\log \log (n)}}{n}} \right\},
\]
then
\[
\Prob^{\, \otimes n}_f\big\{\textsf{R}_\alpha \leqslant 0\big\} \leqslant \beta.
\]
In particular, if $R \in \left[ {\underline R ,\;\overline R } \right]$ with \\[.2cm]
$
\begin{gathered}
\bullet \quad \underline{R} = \left[{\log \log (n)}\right]^{s}\sqrt{\tfrac{\log\log(n)}{n}} \hfill \\
\bullet \quad \overline{R}  = \frac{{n^{2s} }}{{\left[ {\log \log (n)} 
    \right]^{3s + \tfrac{1}{2}} }} \hfill \\ 
\end{gathered} 
$ \\[.2cm]

\noindent then there exists some positive constant $C' = C'(s,\alpha,\beta,M,\|f\|_\infty)$ such that the uniform separation rate of the test $\mathbb{1}_{(0, + \infty )} (\textsf{R}_\alpha)$ over $\mathpzc{A}^s(R,M,G_X)$ satisfies
\[
\rho \left( {\mathbb{1}_{(0, +\infty )}(\textsf{R}_\alpha),\mathpzc{A}^s (R,M,G_X),\beta } \right) \leqslant C' \;R^{\tfrac{1}{{4s + 1}}} \left[ {\frac{{\sqrt {\log \log (n)} }}{n}} \right]^{\tfrac{{2s}}{{2s + 1}}} .
\]
\end{cor}

\begin{remark}
\begin{itemize}
\item The separation rate for the problem of testing ``$f \equiv 0$'' in the classical Gaussian white noise model 
$\mathrm{d}X(t) = f(t)\mathrm{d}t + \epsilon \mathrm{d}W(t)$ has been evaluated for different smoothness classes and distances by Ingster \cite{Ingster:1993a}, Ermakov \cite{Ermakov:1991}, Lepsky and Spokoiny \cite{Lepski:Spoko:1999}, Ingster and Suslina \cite{Ingster:Suslina:1998} (see also the monograph \cite{Ingster:Susl:2002}, and \cite{Lepski:Tsyba:2000} were Lepski and Tsybakov established the asymptotic separation rate -- with constants -- for the $\Lp^\infty$--norm). In \cite{Baraud:2000}, instead, Baraud was able to obtain non--asymptotic bounds on the minimax separation rate in the case of a Gaussian regression model.
From Ingster \cite{Ingster:1993}, we know that the minimax rate of testing over H\"{o}lderian balls $\mathcal{H}^s(R)$ in a Gaussian white noise model is equal to $n^{-2s/(1+4s)}$. From Corollary \ref{cor:CH2:unifSepRates} it seems that we loose a factor equal to $(\log\log(n))^{s/(1+4s)}$ when $s>\tfrac{1}{4}$ but, as Spokoiny proved in \cite{Spokoiny:1996} (see also \cite{Gayraud:Pouet:2005}), adaptivity costs necessarily a logarithmic factor. Therefore we deduce that for $R \in \left[ {\underline R ,\;\overline R } \right]$, our procedure adapts over the approximation space $\mathpzc{A}^s(R,M,G_X)$ at a rate known to be optimal for a particular scale of Besov spaces.
\end{itemize}
\end{remark}

\subsection{Simulation Study}\label{subsec:CH2:SimuStudy}
\begin{figure}[htb]
\begin{center}
\includegraphics[keepaspectratio=true,width=17cm,clip=true]{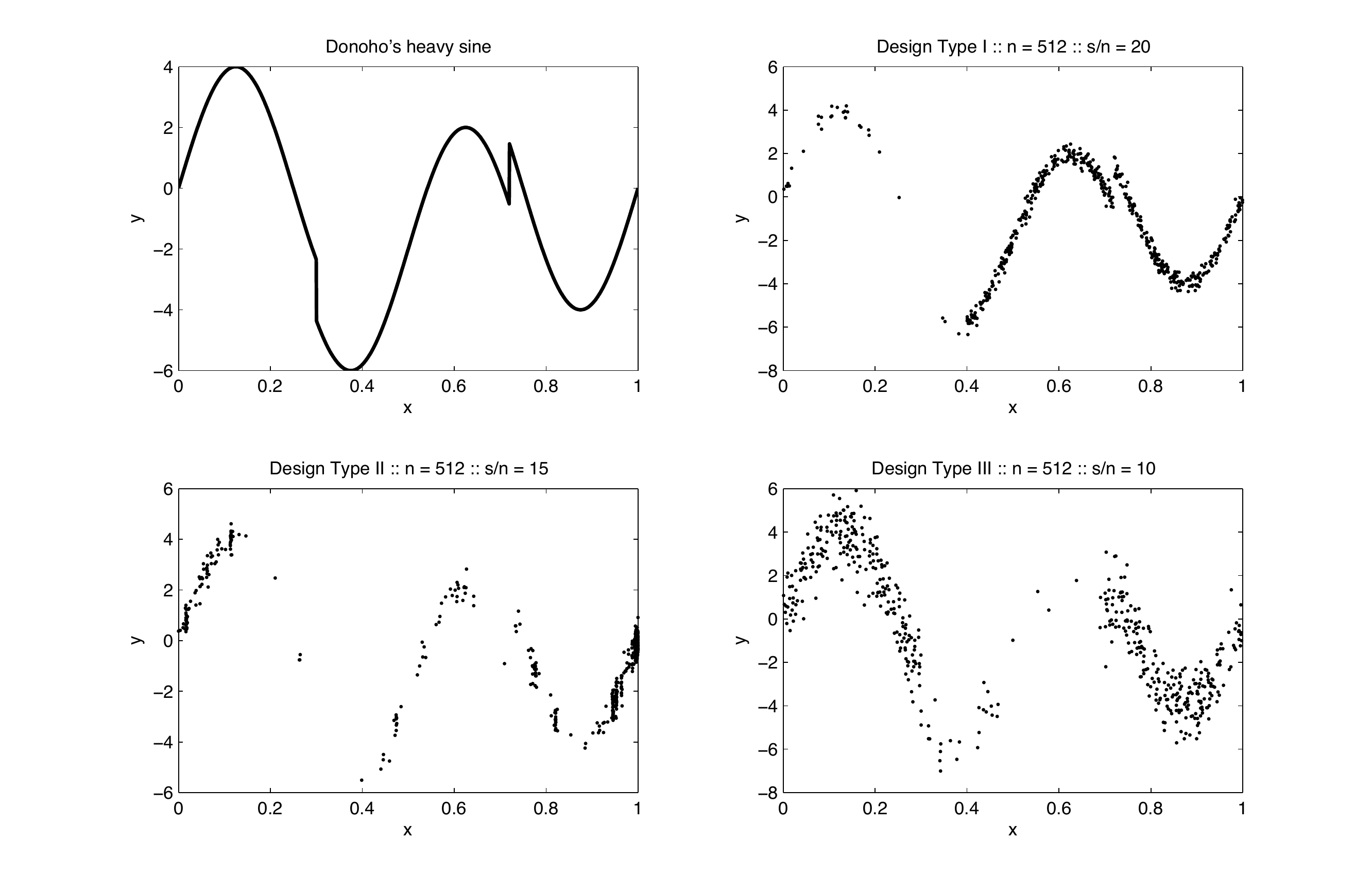}
\end{center}
\caption[Simulation study setup]{\textsf { \small The heavy sine function together with realizations from the three designs chosen to perform the simulation study. With $s/n$ we have denoted the \emph{signal to noise ratio}.
}}
\label{fig:CH2:SimuPic}
\end{figure}

In this section we carry a brief simulation study to evaluate the performances of the proposed testing procedure. Figure \ref{fig:CH2:SimuPic} summarizes the setup. We consider noisy versions of Donoho's \emph{Heavy Sine} function \cite{Dono:John:1994} corresponding to different signal to noise ratios ranging from $10$ to $20$, and three different design distributions that we call \textsf{Type I}, \textsf{II} and \textsf{III}. The sample size is fixed and equal to $512$, whereas we choose to take 
$\textsf{card}(\mathcal{J}_n) = 50$.
Given the nature of the Heavy Sine function, we focus on alternatives of the type
\[
h(x|\kappa) = \kappa \sin(4 \, \pi \,x).
\]
Notice that the true regression function was generated by modifying $h(x|4)$. Finally we set $M = 10$.

As in \cite{Fromont:Laurent:2005}, we have chosen a level $\alpha = 0.05$. The value of $u_\alpha$ and the quantiles $\{r_{n,\MyJ}(u_\alpha)\}_{\MyJ \in \mathcal{J}_n}$ are estimated by $50000$ simulations using a (smoothed) bootstrap procedure. We use $25000$ simulations for the estimation of the $(1-u)$ quantiles $r_{n,\MyJ}(u)$ of the variables
\[
\hat{\textsf{R}}_{n,\MyJ}  = \hat{\theta}_{n,\MyJ}  + \|f_0\|_{\Ldue(G_X)}^2  
    - \frac{2}{n}\sum\nolimits_{i = 1}^n {Y_i f_0 (X_i )},
\]
under the hypothesis ``$f = f_0$'' for $u$ varying on a grid of $(0,\alpha)$, and $25000$ simulations for the estimation of the probabilities
\[
\Prob_{f_0}^{ \otimes n} \left\{ {\mathop {\sup }\limits_{\MyJ \in \mathcal{J}_n } 
    \left[ {\hat \theta _{n,\MyJ} + \|f_0 \|_{\Ldue(G_X )}^2  
    - \frac{2}{n}\sum\limits_i {Y_i f_0 (X_i )}  - r_{n,\MyJ} (u)} \right] > 0} \right\}.
\]
\begin{table}
\centering
\begin{tabular}{ccccccccc} 
\toprule
\textsf{Type I}  \\ 
\midrule
$\kappa = 2$ & $\kappa = 4$ & $\kappa = 6$ & Estim. Lev. \\
\midrule
$0.80$ & $0.77$ & $0.84$ & $0.049$ \\
\toprule
\textsf{Type II}  \\ 
\midrule
$\kappa = 2$ & $\kappa = 4$ & $\kappa = 6$ & Estim. Lev. \\
\midrule
$0.58$ & $0.55$ & $0.60$ & $0.053$ \\
\toprule
\textsf{Type III}  \\ 
\midrule
$\kappa = 2$ & $\kappa = 4$ & $\kappa = 6$ & Estim. Lev. \\
\midrule
$0.44$ & $0.37$ & $0.43$ & $0.052$ \\
\bottomrule
\end{tabular}
\caption[Estimated Power and Level for the test in Section \ref{subsec:CH2:TestUniSepRates}]{\textsf { \small Estimated Power and Level for the test in Section \ref{subsec:CH2:TestUniSepRates}}}\label{tab:CH2:simuRes}
\end{table}
Table \ref{tab:CH2:simuRes} presents the results of our simulation study. 

\section{Discussion}\label{sec:CH2:Discu}

In this short section we collect some remarks regarding the content of this chapter. First of all, it is almost inevitable to mention the most evident weakness of the proposed approach, i.e. the fact that we assumed the design distribution $G_X(\cdot)$ to be completely specified. Although there is a vast literature on the so called \emph{designed experiments} where this type of assumptions are truly welcome, in the present nonparametric regression setting it seems desirable to get rid of it, the most natural way being to assume that $G_X(\cdot)$ belongs to some suitable smoothness class. Clearly this class should necessarily be ``small'' enough so that we are still able to prove the analogs of Theorem \ref{teo:CH2:PowerTest} and Corollary \ref{cor:CH2:unifSepRates}. Notice also that an additional complication we encounter assuming $G_X(\cdot)$ (partially) unknown comes from the fact that now we need to \emph{warp} the initial wavelet basis with some -- possibly smoothed -- version of the empirical distribution function. See the paper \cite{Kerk:Picard:2005} by Picard and Kerkyacharian to have an idea of the intrinsic difficulty of the problem.

Although it is not as disturbing as the previous one, another hypothesis that we might want to relax is the one that requires the knowledge of a (uniform) bound over $|y-f(x)|$. A possible way out here seems to be the use of arguments similar to those adopted by Laurent in \cite{Laurent:2005a} to prove her Proposition 2.

Finally, just a word on the simulation study carried in Section \ref{subsec:CH2:SimuStudy}. Of course this is only a very brief -- although promising -- analysis that can be extended in many directions by considering, for instance, other families of alternatives and regression functions, and possibly a suitable ranges of sample sizes.


\section{Proofs for Section \ref{sec:CH2:GoodOfFit}}\label{sec:CH2:ProofsTest}

\subsection{Proof of Theorem \ref{teo:CH2:PowerTest}}\label{subsec:CH2:ProofTeoPower}

Lets start proving the equality $(\lozenge)$ in Equation (\ref{eq:CH2:riskDecompo}). $\hat{\theta}_{n,k}$ is a U--statistic of order two for $\|\Pi_{\text{S}_k}(f)\|^2_{\Lp^2(G_X)}$ with kernel
\[
h_k(\bz_1,\bz_2) = {\sum\limits_{\ell  \in \mathcal{I}_k } {\left\{ {y_1 e_\ell  (x_1 )} \right\} \cdot \left\{ {y_2 e_\ell (x_2 )} \right\}} },\quad \bz_i = (x_i,y_i),\;i\in\{1,2\}.
\]
Hence,
\[
\hat{\theta}_{n,k}  = \Pi _\emptyset  (h_k) + \frac{2}
{n}\sum\limits_{i = 1}^n {\Pi _{\{ 1\} } (h_k)\left( {\bZ_i } \right)}  + \frac{1}
{{n(n - 1)}}\sum\limits_{i = 2}^n {\sum\limits_{j = 1}^{i - 1} {\Pi _{\{ 1,2\} } (h_k)( {\bZ_i ,\bZ_j })} }, 
\]
where, \\
$
\begin{gathered}
\bfunB
\Pi _\emptyset  (h_k) & = \Exp\big\{ {h_k (\bZ_i ,\bZ_j )} \big\}
    \mathop  = \limits^{\text{indep.}} \sum\limits_{\ell  \in \mathcal{I}_k } 
    {\Exp\big\{ {Y_i e_\ell (X_i)} \big\} \cdot \Exp\big\{ {Y_j e_\ell  (X_j )}\big\}} 
    \mathop = \limits^{\text{id. distr.}} \sum\limits_{\ell  \in \mathcal{I}_k } 
    {\big( {\Exp\left\{ {Ye_\ell  (X)} \right\}} \big)^2 }  = 
    \sum\limits_{\ell  \in \mathcal{I}_k }    {\theta _\ell ^2 } . \hfill \\
\Pi _{\{ 1\} } (h_k)(\bZ_i) & = \Exp\left\{{h_k (\bZ_i,\bZ_j )\big|\bZ_i} \right\} - \Exp\big\{ {h_k (\bZ_i ,\bZ_j )} \big\} = \sum\limits_{\ell  \in \mathcal{I}_k } 
    {\left( {\left\{ {y_i e_\ell  (x_i )} \right\} \cdot \Exp \{ {Y_j e_\ell  (X_j )} \} - \theta _\ell ^2 } \right)} = \hfill \\
{} & = \sum\limits_{\ell  \in \mathcal{I}_k } {\left( {\{ {y_i e_\ell  (x_i )} \} 
    \cdot \theta _\ell   - \theta _\ell ^2 } \right)}  = \sum\limits_{\ell  \in \mathcal{I}_k } {\theta_\ell \big\{ {y_i e_\ell  (x_i ) - \theta _\ell  } \big\}.} \hfill \\
\Pi _{\{1,2\}}(h_k)(\bZ_i,\bZ_j) & = \Exp\big\{ {h_k (\bZ_i ,\bZ_j )\big|\bZ_i ,\bZ_j }\big\} -
    \Exp\big\{ {h_k (\bZ_i ,\bZ_j )\big|\bZ_i } \big\} - 
    \Exp\big\{ {h_k (\bZ_i ,\bZ_j )\big|\bZ_j } \big\} + 
    \Exp\big\{ {h_k (\bZ_i ,\bZ_j )} \big\} =  \hfill \\
{} & = \sum\limits_{\ell \in\mathcal{I}_k } {\left( {\{ {y_i e_\ell  (x_i )} \} \cdot \{ {y_j e_\ell  (x_j )} \} - \theta _\ell  \{ {y_i e_\ell  (x_i )} \} - \theta _\ell \{     {y_j e_\ell (x_j )}\} + \theta _\ell ^2 } \right)} = \hfill \\
{} & = \sum\limits_{\ell  \in \mathcal{I}_k } {\big\{ {y_i e_\ell  (x_i ) - \theta _\ell  } \big\}\cdot \big\{ {y_j e_\ell  (x_j ) - \theta _\ell  } \big\}} .
\efunB
\end{gathered}
$
\\
Hence
\[
\hat{\theta}_{n,k}  = \sum\limits_{\ell  \in \mathcal{I}_k } {\theta _\ell ^2 }  + 
    \frac{2}{n}\sum\limits_i {\sum\limits_{\ell \in \mathcal{I}_k } 
    {\theta _\ell  \big\{ {Y_i e_\ell  (X_i ) - \theta _\ell  } \big\}} } + 
    \frac{1}{{n(n - 1)}}\sum\limits_{i \ne j} {\sum\limits_{\ell  \in \mathcal{I}_k } 
    {\big\{ {Y_i e_\ell (X_i ) - \theta _\ell} \big\}\cdot\big\{ {Y_j e_\ell  (X_j ) - 
    \theta _\ell  } \big\}} }. 
\]
Now note the following equivalences implied by the orthonormality of the system $\{e_i\}_{i\in \mathcal{I}_k}$ in $\Ldue([0,1],G_X)$ \\
$
\begin{gathered}
\bfunB
\| {\Pi _{\text{S}_k } (f)} \|_{\Ldue(G_X )}^2 & = 
    \left\|{\sum\limits_{\ell  \in \mathcal{I}_k }{\theta _\ell \,e_\ell}} \right\|_{\Ldue(G_X)} 
    \mathop  = \limits^{\text{orthonorm.}} 
    \sum\limits_{\ell  \in \mathcal{I}_k } {\theta _\ell ^2 }. \hfill \\
2 \left( {\EmpProb_n  - \Prob} \right)\big(\Pi _{\text{S}_k } (f)\big) & = 
    2 \left\{ {\frac{1}{n}\sum\limits_{i = 1}^n {Y_i \,\Pi _{\text{S}_k } (f)(X_i )}  - 
    \left\langle {\Pi _{\text{S}_k } (f),f} \right\rangle _{\Ldue(G_X )} } \right\} =  \hfill \\
& = 2\left\{ {\frac{1}{n}\sum\limits_{i = 1}^n {Y_i 
    \left[ {\sum\limits_{\ell  \in \mathcal{I}_k } {\theta _\ell  \,e_\ell (X_i)} } \right]}  
    - \left\| {\Pi _{\text{S}_k } (f)} \right\|_{\Ldue(G_X )}^2 } \right\} = \hfill \\
& = 2 \left\{ {\frac{1}{n}\sum\limits_{i = 1}^n {\left[ {\sum\limits_{\ell  \in \mathcal{I}_k }
    {\theta _\ell  \,Y_i \,e_\ell (X_i )} } \right]}  - \sum\limits_{\ell  \in \mathcal{I}_k } 
    {\theta _\ell ^2 } } \right\} =  \hfill \\
& = 2\left\{ {\frac{1}{n}\sum\limits_{i = 1}^n {\left[ {\sum\limits_{\ell  \in \mathcal{I}_k } {\theta _\ell  \,Y_i \,e_\ell (X_i )} } \right]}  - \frac{1}{n}\sum\limits_{i = 1}^n 
    {\left[ {\sum\limits_{\ell  \in \mathcal{I}_k } {\theta _\ell ^2 } } \right]} } \right\} = \hfill \\
& = \frac{2}{n}\sum\limits_{i = 1}^n {\sum\limits_{\ell \in \mathcal{I}_k } {\theta _\ell \,
    \big\{ Y_i \,e_\ell  (X_i ) - \theta _\ell  \big\} } }.
\efunB
\end{gathered} 
$ \\
So 
\[
\hat{\theta}_{n,k}  = \left\| {\Pi _{\text{S}_k } (f)} \right\|_{\Ldue(G_X )}^2  + 2\left( {\EmpProb_n - \Prob} \right)\big(\Pi _{\text{S}_k } (f)\big) + \widetilde{U}_{n,k}. 
\]
Finally, by using the fact that
\begin{eqnarray*}
\|f - \Pi_{\text{S}_k}(f)\|^2_{\Ldue(G_X)} &=& \|f\|^2_{\Ldue(G_X)} + \|\Pi_{\text{S}_k}(f)\|^2_{\Ldue(G_X)} - 2\langle f, \Pi_{\text{S}_k}(f) 
    \rangle_{\Ldue(G_X )}^2 = \\
&=& \|f\|^2_{\Ldue(G_X)} + \|\Pi_{\text{S}_k}(f)\|^2_{\Ldue(G_X)} - 2 \|\Pi_{\text{S}_k}(f)\|^2_{\Ldue(G_X)} = \\
&=&  \|f\|^2_{\Ldue(G_X)} - \|\Pi_{\text{S}_k}(f)\|^2_{\Ldue(G_X)},
\end{eqnarray*}
from Equation (\ref{eq:CH2:riskDecompo}) we obtain
\begin{eqnarray*}
\hat{\textsf{R}}_{n,k}  &=& \hat{\theta}_{n,k}  + \| {f_0 } \|_{\Ldue(G_X )}^2  - \frac{2}{n}\sum\limits_i {Y_i f_0 (X_i )}  =  \hfill \\
&=& \left\{ {\widetilde{U}_{n,k}  + 2\left( {\EmpProb_n  - \Prob} \right)
    \big(\Pi_{\text{S}_k }(f) \pm f\big) - \big[ { \pm \| f \|_{\Ldue(G_X )}^2  
    - \| {\Pi _{\text{S}_k } (f)} \|_{\Ldue(G_X )}^2 } \big]} \right\} + 
    \| {f_0 } \|_{\Ldue(G_X )}^2  - 2\EmpProb_n (f_0 ) =  \hfill \\
&=& \left\{ {\widetilde{U}_{n,k}  + 2\left( {\EmpProb_n  - \Prob} \right)
   \big(\Pi _{\text{S}_k }(f) - f\big) - \| {f -\Pi _{\text{S}_k } (f)} \|_{\Ldue(G_X )}^2} 
   \right\} + \hfill \\
& & + \| {f_0 } \|_{\Ldue(G_X )}^2  + \| f \|_{\Ldue(G_X )}^2  \pm 2 \Prob(f) + 2(\EmpProb_n - \Prob)(f) - 2\EmpProb_n (f_0 ) =  \hfill \\
&=& \big\{ {\; \cdot \;} \big\} + \big\{ {\| {f_0 } \|_{\Ldue(G_X )}^2  + 
    \| f \|_{\Ldue(G_X )}^2  - 2\Prob(f)} \big\} + 
    \big\{ {2(\EmpProb_n  - \Prob)(f) - 2(\EmpProb_n  - \Prob)(f_0 )} \big\} =  \hfill \\
&=& \widetilde{U}_{n,k}  + 2\left( {\EmpProb_n  - \Prob} \right)\big(\Pi _{\text{S}_k } (f) - f\big) - \| {f - \Pi _{\text{S}_k } (f)} \|_{\Ldue(G_X )}^2  
    + 2(\EmpProb_n  - \Prob)\big(f - f_0\big) + \| {f - f_0 } \|_{\Ldue(G_X )}^2,
\end{eqnarray*}

\noindent and this complete the proof.

Now, given $\beta \in (0,1)$  we know that
\[
\Prob_f^{ \otimes n} \left\{ {\textsf{R}_\alpha   \leqslant 0} \right\} = 
    \Prob_f^{ \otimes n}\left\{ {\mathop {\sup }\limits_{\MyJ \in \mathcal{J}_n } 
    \left[ {\hat{\theta}_{n,\MyJ} + \|f_0 \|_{\Ldue(G_X )}^2  - 
    \tfrac{2}{n}\sum\nolimits_i {Y_i f_0 (X_i )}  - r_{n,\MyJ} (u_\alpha  )} \right] 
    \leqslant 0} \right\},
\]
hence
\begin{eqnarray}\label{eq:CH2:ProxyPower}
\Prob_f^{ \otimes n} \left\{ {\textsf{R}_\alpha   \leqslant 0} \right\} 
    &\leqslant& \mathop {\inf }\limits_{\MyJ \in \mathcal{J}_n } 
    \Prob_f^{ \otimes n} \left\{ {\hat{\theta}_{n,\MyJ}  + \|f_0 \|_{\Ldue(G_X )}^2     - \tfrac{2}{n}\sum\limits_i {Y_i f_0 (X_i )}  - r_{n,\MyJ} (u_\alpha) \leqslant 0} \right\} = \nonumber \hfill \\
&=& \mathop {\inf }\limits_{\MyJ \in \mathcal{J}_n } 
    \Prob_f^{ \otimes n} \left\{ {\widetilde{U}_{n,\MyJ} + 
    2\left( {\EmpProb_n  - \Prob} \right)\big(\Pi _{\text{S}_\MyJ } (f) - f\big) - 
    \| {f - \Pi _{\text{S}_\MyJ} (f)}\|_{\Ldue(G_X )}^2 + } \right. \hfill \\
& & \quad\quad\quad\quad + \left. {2(\EmpProb_n  - \Prob)(f - f_0 ) + \| {f - f_0 } \|_{\Ldue(G_X )}^2 - r_{n,\MyJ} (u_\alpha  ) \leqslant 0} \right\}. \nonumber
\end{eqnarray}
Following \cite{Fromont:Laurent:2005}, we will split the control of the power in three steps, involving separately $\widetilde{U}_{n,\MyJ}$, $2(\EmpProb_n - \Prob)\big(\Pi_{\text{S}_\MyJ}(f) - f \big)$, and $2(\EmpProb_n - \Prob)\big(f-f_0\big)$. To handle the last two terms, we will use the following version of the Bernstein's inequality with constants provided by Birg\'{e} and Massart in \cite{Birge:Massart:1998}:
\begin{lem}\label{lem:CH2:BernesteinIneq}
Let $\{U_i\}_{i\in\nSeq}$ be independent random variables such that for all $i\in\nSeq$
\begin{itemize}
\item $|U_i| \leqslant b$,
\item $\Exp(U_i^2) \leqslant \delta^2$.
\end{itemize}
Then, for all $u>0$,
\bequ\label{eq:CH2:BernsteinIneq}
\Prob^{ \otimes n} \left\{ {\frac{1}{n}\sum\limits_{i = 1}^n 
{\left[ {U_i  - \Exp(U_i )} \right]} > \frac{{\delta \sqrt {2u} }}{{\sqrt n }} 
+ \frac{{b\,u}}{{3n}}} \right\} \leqslant \mathrm{e}^{ - u}.
\eequ
\end{lem}

\subsubsection{$\bullet \quad $ Control of $\widetilde{U}_{n,\MyJ}$}\label{subsubsec:CH2:ControlUstat}
We start with the following lemma whose proof is postponed to Section \ref{subsec:CH2:proofHoudre}.
\begin{lem}\label{lem:CH2:HoudreReyn_inAction}
Under the hypotheses and using the notation of Theorem \ref{teo:CH2:PowerTest}, there exists an absolute constant $\kappa_0$ such that, for all $u > 0$, we have
\bequ\label{eq:CH2:Houdre_inAction}
\Prob_f^{ \otimes n} \left\{ {\left| {\widetilde{U}_{n,\MyJ} } \right| > 
\frac{{\kappa _0 }}{n}\left[ {\tau _\infty  \sqrt {u 2^{\MyJ} }  + \tau _\infty  u + M^{\,2} \frac{{u^2 2^{\MyJ} }}{n}} \right]} \right\} \leqslant 5.6\,\mathrm{e}^{ - u} .
\eequ
\end{lem}
Now, let  $u_{\textsf{I}} = u_{\textsf{I}}(\beta) = \log(3/\beta)$, and $u_{\textsf{II}} = u_{\textsf{II}}(\beta) = u_{\textsf{I}} + \log(5.6)$. Then, from Equation (\ref{eq:CH2:Houdre_inAction}), we obtain
\bequ\label{eq:CH2:Houdre_inAction_II}
\Prob_f^{ \otimes n} \left\{ {\left| {\widetilde{U}_{n,\MyJ} } \right| <  
    - \frac{{\kappa _0 }}{n}\left[ {\tau _\infty  
    \sqrt {u_{\textsf{II}} 2^{\MyJ} }  + \tau _\infty  u_{\textsf{II}}  
    + M^{\,2} \frac{{u_{\textsf{II}}^2 2^{\MyJ} }}{n}} \right]} \right\} 
    \leqslant \frac{\beta }{3}.
\eequ
where $\tau_\infty = \|f\|^2_\infty + \|\sigma^2\|_\infty$.

\subsubsection{$\bullet \quad $Control of $2(\EmpProb_n - \Prob)\big(\Pi_{\text{S}_\MyJ}(f) - f \big)$}\label{subsubsec:CH2:ControlPartA}

In order to apply Lemma \ref{lem:CH2:BernesteinIneq}, let
\[
U = 2Y\,\big[\Pi_{\text{S}_\MyJ}(f)(X) - f(X)\big],
\]
then
\begin{itemize}
\item $|U| = \left| {2Y\left[ {\Pi _{S_\MyJ } (f)(X) - f(X)} \right]} \right| \leqslant 2M\left\{ {\mathop {\sup }\limits_{x \in [0,1]} \left| {\Pi _{S_{\MyJ} } (f)(x)} \right| - \mathop {\sup }\limits_{x \in [0,1]} \left| {f(x)} \right|} \right\} = 4M\,\|f\|_\infty.$
\item 
$\bfunB
    \Exp\left( {U^2 } \right) = &
    \Exp\left\{{4Y^2 \left({\Pi_{\text{S}_{\MyJ}}(f) - f} \right)^2(X)}\right\} = \hfill \\
{ } & = 4\Exp\left\{ {\left[ {f^2 (X) + \sigma ^2 (X)} \right]
    \left( {\Pi _{\text{S}_{\MyJ} }(f) - f} \right)^2 (X)} \right\} 
    \leqslant 4\tau_\infty\Exp\left\{{\left({\Pi_{\text{S}_{\MyJ}} (f) - f}\right)^2(X)}
    \right\}  =  \hfill \\
{ } & = 4\tau _\infty  \|\Pi _{\text{S}_{\MyJ} } (f) - f\|_{\Ldue(G_X )}^2 . \hfill \\ 
\efunB
$
\end{itemize}
Hence, applying Lemma \ref{lem:CH2:BernesteinIneq}, we have
\[
\Prob_f^{ \otimes n} \left\{ {2( {\EmpProb_n  - \Prob} )
    \big( {\Pi _{\text{S}_{\MyJ} } (f) - f} \big) < 
    - \frac{{2\sqrt {2\tau _\infty  } \sqrt u }}{{\sqrt n }} \,
    \|\Pi _{\text{S}_{\MyJ} } (f) - f\|_{\Ldue(G_X )}  - 
    \frac{{4M \, \|f\|_\infty  u}}{{3n}}} \right\} \leqslant \mathrm{e}^{ - u}. 
\]
By the inequality $2 a b \leqslant \tfrac{4}{\gamma}a^2 + \tfrac{\gamma}{4} b^2$ we then have
\[
2\left[ {\sqrt {\frac{{2\tau _\infty  u}}{n}} } \right] 
    \cdot \left[ {\|\Pi _{\text{S}_{\MyJ} } (f) - f\|_{\Ldue(G_X )} } \right] \leqslant 
    \frac{4}{\gamma }\frac{{2\tau _\infty  u}}{n} + \frac{\gamma }{4}
    \|\Pi _{\text{S}_{\MyJ} } (f) - f\|^2_{\Ldue(G_X )} ,
\]
and consequently
\[
\Prob_f^{ \otimes n} \left\{ {2\left( {\EmpProb_n  - \Prob} \right)
    \left( {\Pi _{\text{S}_{\MyJ} } (f) - f} \right) 
    + \frac{\gamma }{4}\|\Pi _{\text{S}_{\MyJ} } (f) - f\|_{\Ldue(G_X )}^2  <  
    - \left[ {\frac{8}{\gamma }\tau _\infty   
    + \frac{4}{3}M\, \|f\|_\infty  } \right]\frac{u}{n}} \right\} \leqslant \mathrm{e}^{ - u}. 
\]
Finally, taking $u_{\textsf{I}} = u_{\textsf{I}}(\beta) = \log(3/\beta)$, as before, we get
\bequ\label{eq:CH2:Bernstein_inAction_I}
\Prob_f^{ \otimes n} \left\{ {2\left( {\EmpProb_n  - \Prob} \right)
    \left( {\Pi _{\text{S}_{\MyJ} } (f) - f} \right) 
    + \frac{\gamma }{4}\|\Pi _{\text{S}_{\MyJ} } (f) - f\|_{\Ldue(G_X )}^2  <  
    - \left[ {\frac{8}{\gamma }\tau _\infty   
    + \frac{4}{3}M \, \|f\|_\infty } \right]\frac{u_{\textsf{I}}}{n}} \right\} 
    \leqslant \frac{\beta}{3}.
\eequ

\subsubsection{$\bullet \quad $Control of $2(\EmpProb_n - \Prob)\big(f-f_0\big)$}\label{subsubsec:CH2:ControlPartB}
Proceeding as in the previous section, let
\[
U = 2Y\,\big[f(X) - f_0(X)\big],
\]
then
\begin{itemize}
\item $|U| = \left| {2Y\left[ f(X) - f_0(X) \right]} \right| 
    \leqslant 2M\big\{ \|f\|_\infty + \|f_0\|_\infty \big\}$,
\item 
$\bfunB
    \Exp\left( {U^2 } \right) = &
    \Exp\left\{ {4Y^2 \left( {f - f_0} \right)^2 (X)} \right\} =  \hfill \\
{ } & = 4\Exp\left\{ {\left[ {f^2 (X) + \sigma ^2 (X)} \right]
    \left( {f - f_0} \right)^2 (X)} \right\} 
    \leqslant 4\tau _\infty \Exp \left\{ {\left( {f - f_0} \right)^2(X)} \right\}  =  \hfill \\
{ } & = 4\tau _\infty  \|f - f_0\|_{\Ldue(G_X )}^2 . \hfill \\ 
\efunB
$
\end{itemize}
Hence, applying Lemma \ref{lem:CH2:BernesteinIneq}, we have
\[
\Prob_f^{ \otimes n} \left\{ {2( {\EmpProb_n  - \Prob} )
    \big( {f - f_0} \big) < - \frac{{2\sqrt {2\tau _\infty  } \sqrt u }}{{\sqrt n }} \,
    \|f - f_0 \|_{\Ldue(G_X )}  - 
    \frac{{2M \, \{\|f\|_\infty + \|f_0\|_\infty\} u}}{{3n}}} \right\} \leqslant \mathrm{e}^{ - u}. 
\]
Applying again the inequality $2 a b \leqslant \tfrac{4}{\gamma}a^2 + \tfrac{\gamma}{4} b^2$ we then have
\[
2\left[ {\sqrt {\frac{{2\tau _\infty  u}}{n}} } \right] 
    \cdot \left[ {\|f - f_0\|_{\Ldue(G_X )} } \right] \leqslant 
    \frac{4}{\gamma }\, \frac{{2\tau _\infty  u}}{n} + \frac{\gamma }{4} \,
    \|f - f_0\|^2_{\Ldue(G_X )} ,
\]
and consequently
\[
\Prob_f^{ \otimes n} \left\{ {2\left( {\EmpProb_n  - \Prob} \right)
    \left( {f - f_0} \right) 
    + \frac{\gamma }{4}\|f - f_0\|_{\Ldue(G_X )}^2  <  
    - \left[ {\frac{8}{\gamma }\tau _\infty   
    + \frac{2}{3}M \big\{\|f\|_\infty + \|f_0\|_\infty\big\}} \right]\frac{u}{n}} \right\} 
    \leqslant \mathrm{e}^{ - u}. 
\]
Finally, taking $u_{\textsf{I}} = u_{\textsf{I}}(\beta) = \log(3/\beta)$, as always, we get
\bequ\label{eq:CH2:Bernstein_inAction_II}
\Prob_f^{ \otimes n} \left\{ {2\left( {\EmpProb_n  - \Prob} \right)
    \left( {f - f_0} \right) 
    + \frac{\gamma }{4}\|f - f_0\|_{\Ldue(G_X )}^2  <  
    - \left[ {\frac{8}{\gamma }\tau _\infty   
    + \frac{2}{3}M \big\{\|f\|_\infty + \|f_0\|_\infty\big\}} \right]
    \frac{u_{\textsf{I}}}{n}} \right\} \leqslant \frac{\beta}{3}.
\eequ

\subsubsection{$\bullet \quad $ Conclusion}\label{subsubsec:CH2:PowerConcludeProof}

Combining Equation (\ref{eq:CH2:ProxyPower}) with the bounds presented in Equations (\ref{eq:CH2:Houdre_inAction_II},\ref{eq:CH2:Bernstein_inAction_I},\ref{eq:CH2:Bernstein_inAction_II}),
we get
\begin{eqnarray*}
\Prob_f^{ \otimes n} \big\{ {\textsf{R}_\alpha \leqslant 0} \big\} \leqslant \beta &+&
    \mathop {\inf }\limits_{\MyJ \in \mathcal{J}_n }\Prob_f^{ \otimes n}
    \Big\{ {\|f - f_0 \|_{\Ldue(G_X )}^2 \leqslant
    \|\Pi _{\text{S}_{\MyJ} } (f) - f\|_{\Ldue(G_X )}^2  
    + r_{n,\MyJ} (u_\alpha  ) + } \big. \hfill \\
& & \quad\quad\quad\quad + \tfrac{{\kappa _0 }}{n}
    \big[\tau _\infty  
    \sqrt {u_{\textsf{II}} 2^{\MyJ} } + \tau _\infty u_{\textsf{II}} 
    + M^{\,2} u_{\textsf{II}}^2 2^{\MyJ}\tfrac{1}{n}\big]+\hfill \\
& & \quad\quad\quad\quad + \tfrac{\gamma }{4} 
    \|\Pi _{\text{S}_{\MyJ} } (f) - f\|_{\Ldue(G_X )}^2 + 
    \big[\tfrac{8}{\gamma }\tau _\infty   + \tfrac{4}{3}M\,\|f\|_\infty  \big]
    \tfrac{{u_\textsf{I} }}{n} +  \hfill \\
& & \quad\quad\quad\quad + \left. { \tfrac{\gamma }{4}\|f - f_0 \|_{\Ldue(G_X )}^2  
  + \big[\tfrac{8}{\gamma }\tau _\infty   + \tfrac{2}{3}M\,
  \{ \|f\|_\infty   + \|f_0\|_\infty  \}\big]\tfrac{{u_\textsf{I} }}{n}} \right\}.
\end{eqnarray*}
So, if exists $\MyJ \in \mathcal{J}_n$ such that
\begin{eqnarray*}
(1 - \tfrac{\gamma }{4}) \, \|f - f_0 \|_{\Ldue(G_X )}^2  &>& 
    (1 + \tfrac{\gamma }{4})\, \|\Pi _{\text{S}_{\MyJ} } (f) - f\|_{\Ldue(G_X )}^2  
    + \tfrac{{\kappa _0 }}{n}\left[ {\tau _\infty  \sqrt {u_{\textsf{II}} 2^{\MyJ} }  + 
    \tau _\infty  u_{\textsf{II}}  + M^{\,2} u_{\textsf{II}}^2 2^{\MyJ} 
    \tfrac{1}{n}} \right] +  \hfill \\
&+& \left[ \tfrac{16}{\gamma} \tau _\infty   + 
    2 M \big( \|f\|_\infty + \tfrac{1}{3} \|f_0\|_\infty \big) \right]
    \tfrac{{u_\textsf{I}}}{n} + r_{n,\MyJ} (u_\alpha),
\end{eqnarray*}
then 
\[
\Prob_f^{ \otimes n} \big\{ {\textsf{R}_\alpha   \leqslant 0} \big\} \leqslant \beta,
\]
and this complete the proof of Theorem \ref{teo:CH2:PowerTest}.

\subsection{Proof of Corollary \ref{cor:CH2:unifSepRates}}\label{subsec:CH2:ProofSepRates}

We will split the proof of Corollary \ref{cor:CH2:unifSepRates} in two parts: in the first one we will bound $r_{n,\MyJ}(u_\alpha)$, the $1-u_\alpha$ quantile under the null hypothesis of the test statistic $\hat{\textsf{R}}_{n,\MyJ}$; whereas in the second one, we shall use this bound together with Theorem \ref{teo:CH2:PowerTest} to complete the proof.

\subsubsection{$\bullet \quad $ Upper bound for $r_{n,\MyJ}(u_\alpha)$, $\MyJ \in \mathcal{J}_n$}\label{subsubsec:CH2:QuantileUpperBound}

In this section we will prove the following lemma:
\begin{lem}\label{lem:CH2:boundTheQuantile}
Under the hypotheses and using the notation of Corollary \ref{cor:CH2:unifSepRates}, there exists a positive constant $C(\alpha)$, such that
\[
r_{n,\MyJ}(u_\alpha) \leqslant \widetilde{r}_{n,\MyJ}(\alpha),
\]
where
\[
\widetilde{r}_{n,\MyJ} (\alpha ) = \frac{{C(\alpha )}}{n}
\left\{ 
    {\tau _{0,\infty } 2^{\MyJ/2} \sqrt {\log \log (n)}  
   + {2\left[ {\tau _{0,\infty }  + \tfrac{1} {3}M \|f_0\|_\infty  } \right]}  \log \log (n) 
   + M^{\,2} 2^{\MyJ} \frac{{\left[ {\log \log (n)} \right]^2 }}{n}} \right\},
\]
\end{lem}
with $\tau_{0,\infty} = \|f_0\|^2_\infty + \|\sigma^2\|_\infty$.
\begin{proof}
First of all notice that, by hypothesis, and for all $n\in\MyNN$,
\[
\mathcal{J}_n = \big\{0,\ldots,\log _2 \left\{ {\tfrac{{n^2 }}{{[\log \log (n)]^3 }}} \right\}\big\} \Rightarrow 
\textsf{card}(\mathcal{J}_n) = 1 + \log_2\big\{ \tfrac{n^2}{[\log\log(n)]^3}\big\} \leqslant
1 + \log_2(n^2).
\]
Hence, under the null ``$f = f_0$'', and for $\alpha_n = \alpha/[1+\log_2(n^2)]$, we get
\[
\Prob_{f_0 }^{ \otimes n} \left\{ {\mathop {\sup }\limits_{\MyJ \in \mathcal{J}_n } 
    \big\{ \hat{\textsf{R}}_{n,\MyJ} - r_{n,\MyJ}(\alpha _n )\big\}  > 0} \right\} 
    \leqslant \sum\limits_{\MyJ \in \mathcal{J}_n } {\Prob_{f_0 }^{ \otimes n} 
    \left\{ {\hat{\textsf{R}}_{n,\MyJ}  - r_{n,\MyJ} (\alpha _n ) > 0} \right\}} 
    \leqslant \sum\limits_{\MyJ \in \mathcal{J}_n } {\frac{\alpha }{{[1 + \log _2 (n^2 )]}}}  
    \leqslant \alpha .
\]
Consequently, 
\[
\alpha _n  \leqslant \sup \left\{ {u \in (0,1):\Prob_{f_0 }^{ \otimes n} 
    \left[ {\mathop {\sup }\limits_{\MyJ \in \mathcal{J}_n } 
    \left\{ {\hat{\textsf{R}}_{n,\MyJ}  - r_{n,\MyJ} (u)} \right\} > 0} \right] 
    \leqslant \alpha } \right\} = u_\alpha \Rightarrow r_{n,\MyJ} (u_\alpha  ) 
    \leqslant r_{n,\MyJ} (\alpha _n ).
\]
Hence, all we have to do is to find an upper bound for $r_{n,\MyJ} (\alpha _n )$.

Working under the null, from Equation (\ref{eq:CH2:riskDecompo}) we obtain
\[
\hat{\textsf{R}}_{n,\MyJ}  = \widetilde{U}_{n,\MyJ}  
    + 2\left( {\EmpProb_n  - \Prob} \right)\big(\Pi _{\text{S}_\MyJ } (f_0 ) - f_0 \big) 
    - \|\Pi _{\text{S}_\MyJ } (f_0 ) - f_0\|_{\Ldue(G_X )}^2.
\]
At this point, once we set $u_{n,\textsf{I}} = \log(2/\alpha_n)$ and $u_{n,\textsf{II}} = u_{n,\textsf{I}} + \log(5.6)$, we can proceed as in the proof of Theorem \ref{teo:CH2:PowerTest} obtaining the following bounds
\begin{itemize}
\item By Lemma \ref{lem:CH2:HoudreReyn_inAction},
\[
\Prob_{f_0 }^{ \otimes n} \left\{ {\left| {\widetilde{U}_{n,\MyJ} } \right| > 
    \frac{{\kappa _0 }}{n}\left[ {\tau _{0,\infty } \sqrt {u_{n,\textsf{II}} 2^{\MyJ} }  
    + \tau _{0,\infty } u_{n,\textsf{II}}  + M^{\,2} \frac{{u_{n,\textsf{II}}^2 2^{\MyJ} }}{n}}
    \right]} \right\} \leqslant \frac{{\alpha _n }}{2},
\]
\item By Lemma \ref{lem:CH2:BernesteinIneq}, and using the inequality $2 a b \leqslant a^2 + b^2$,
\[
\Prob_{f_0 }^{ \otimes n} \left\{ {2\left( {\EmpProb_n  - \Prob} \right)
    \left( {\Pi_{\text{S}_{\MyJ} } (f_0 ) - f_0 } \right) 
    - \|\Pi _{\text{S}_{\MyJ} } (f_0 ) - f_0 \|_{\Ldue(G_X )}^2 > 
    2\left[ {\tau _{0,\infty} + \frac{2}{3}M\,\|f_0\|_\infty  } \right]
    \frac{{u_{n,\textsf{I}} }}{n}} \right\} \leqslant \frac{{\alpha _n }}{2}.
\]
\end{itemize}
Combining these two inequalities we get
\[
\Prob_{f_0 }^{ \otimes n} \left\{ {\hat{\textsf{R}}_{n,\MyJ}  > 
    \tfrac{{\kappa _0 }}{n}\left[ {\tau _{0,\infty } \sqrt {u_{n,\textsf{\textsf{II}}} 2^{\MyJ} }  
    + \tau _{0,\infty } u_{n,\textsf{II}}  + M^{\,2} u_{n,\textsf{II}}^2 2^{\MyJ} \tfrac{1}{n}} \right] 
    + 2\left[ {\tau _{0,\infty} + \tfrac{2}{3}M\,\|f_0\|_\infty  } \right]u_{n,\textsf{I}} 
    \tfrac{1}{n}} \right\} \leqslant \alpha _n .
\]
Finally, it is easy to see that we can find two constants $C'(\alpha)$ and $C''(\alpha)$ such that $u_{n,\textsf{I}} \leqslant C'(\alpha) \log\log(n)$ and $u_{n,\textsf{II}} \leqslant C''(\alpha) \log\log(n)$, therefore
{\scriptsize{
\begin{eqnarray*}
& & \tfrac{{\kappa _0 }}{n}\left[ {\tau _{0,\infty } \sqrt {u_{n,\textsf{II}} 2^{\MyJ} }  
    + \tau _{0,\infty } u_{n,\textsf{II}}  + M^{\,2} 2^{\MyJ} \tfrac{{u_{n,\textsf{II}}^2 }}{n}}
    \right] + 2\left[ {\tau _{0,\infty }  + \tfrac{2}{3} M \|f_0\|_\infty  } \right]
    \tfrac{{u_{n,\textsf{I}} }}{n} \leqslant  \hfill \\
&\leqslant& \tfrac{1}{n}\left\{ {\left\{ {\kappa _0 \sqrt {C''(\alpha )} } \right\}
    \tau _{0,\infty } 2^{\MyJ/2} \sqrt {\log \log (n)}  + \left\{ {\kappa _0 
    C''(\alpha )}\right\} \tau _{0,\infty } \log \log (n) + \left\{ {\kappa _0 
    [C''(\alpha )]^2 } \right\}M^{\,2} 2^{\MyJ} \tfrac{{[\log \log (n)]^2 }}{n} 
    + \left\{ {2C'(\alpha )} \right\}\left[ {\tau _{0,\infty }  + \tfrac{2}{3} M \|f_0\|_\infty}
    \right]\log \log (n)} \right\} \leqslant  \hfill \\
&\leqslant& \tfrac{{C(\alpha )}}{n}\left\{ {\tau _{0,\infty } 2^{\MyJ/2} \sqrt {\log \log (n)}
    + \tau _{0,\infty } \log \log (n) + M^{\,2} 2^{\MyJ} \tfrac{{[\log \log (n)]^2 }}{n} 
    + \left[ {\tau _{0,\infty }  + \tfrac{2}{3} M \|f_0\|_\infty  } \right]
    \log \log (n)} \right\} =  \hfill \\
&=& \tfrac{{C(\alpha )}}{n}\left\{ {\tau _{0,\infty } 2^{\MyJ/2} \sqrt {\log \log (n)}  
    + 2\left[ {\tau _{0,\infty }  + \tfrac{1}{3} M \|f_0\|_\infty  } \right]
    \log \log (n) + M^{\,2} 2^{\MyJ} \tfrac{{[\log \log (n)]^2 }}{n}} \right\},
\end{eqnarray*}
}}
where $C(\alpha ) = \max \left\{ {\kappa _0 \sqrt {C''(\alpha )} ,\kappa _0 C''(\alpha ),\kappa _0 [C''(\alpha )]^2 ,2C'(\alpha )} \right\}$. And this complete the proof.
\end{proof}

\subsubsection{$\bullet \quad $ Separation rates}\label{subsubsec:CH2:SepRates}

Combining Theorem \ref{teo:CH2:PowerTest} and Lemma \ref{lem:CH2:boundTheQuantile}, for each $\beta \in (0,1)$ we get that
\[
\Prob_f^{ \otimes n} \left\{ {\textsf{R}_\alpha   \leqslant 0} \right\} \leqslant \beta, 
\]
for every $f(\cdot)$ such that
\[
\|f - f_0\|_{\Ldue(G_X )}^2  > (1 + \gamma )\;\mathop {\inf }\limits_{\MyJ \in \mathcal{J}_n } \left\{ {\|f - \Pi _{\text{S}_\MyJ } (f)|_{\Ldue(G_X )}^2  + 
    \widetilde{r}_{n,\MyJ} (\alpha ) + \textsf{V}_{n,\MyJ} (\beta )} \right\}.
\]
Now, assuming that $f \in \mathpzc{A}^2(R,M,G_X)$, the right hand side in the last equation reduces to
\begin{eqnarray*}
& & \mathop {\inf }\limits_{\MyJ \in \mathcal{J}_n } \left\{ {R^2 2^{ - 2\MyJ\,s}  
    + \tfrac{{C(\alpha )}}{n}\tau _{0,\infty } 2^{\MyJ/2} \sqrt {\log \log (n)}  
    + \tfrac{{2C(\alpha )}}{n}\left[ {\tau _{0,\infty } + \tfrac{1}{3}M\|f_0 \|_\infty}\right]
    \log \log (n) + \tfrac{{C(\alpha )}}{n}M^{\,2} 2^{\MyJ} \tfrac{{[\log \log (n)]^2 }}{n} + }
    \right. \hfill \\
& & \quad\quad\quad
    \biggl. { + \tfrac{{C_1 }}{n}[\tau _\infty  2^{\MyJ/2}  + \tfrac{{M^2 }}{n}2^{\MyJ} ] 
    + \tfrac{{C_2 }}{n}} \biggr\} =  \hfill \\
&=& \mathop {\inf }\limits_{\MyJ \in \mathcal{J}_n } \left\{ {R^2 2^{ - 2\MyJ\,s} 
    + C_1 \tau _\infty \tfrac{{\sqrt {2^{\MyJ}  \cdot 1} }}{n} 
    + C(\alpha )\tau _{0,\infty } \tfrac{{\sqrt {2^{\MyJ} \log \log (n)} }}{n} + 
    C_1 M^{\,2} \tfrac{{2^{\MyJ}  \cdot 1}}{{n^2 }} + C(\alpha )M^{\,2} \tfrac{{2^{\MyJ} 
    [\log \log (n)]^2 }}{{n^2 }} + } \right. \hfill \\
& & \quad\quad\quad
    \Biggl. { + 2C(\alpha )\left[ {\tau _{0,\infty }  + \tfrac{1}{3}M\|f_0\|_\infty} \right]
    \tfrac{{\log \log (n)}}{n} + \tfrac{{C_2  \cdot 1}}{n}} \Biggr\} \leqslant  \hfill \\
&\mathop \leqslant \limits^{(\diamondsuit )}& \mathop {\inf }\limits_{\MyJ \in \mathcal{J}_n } \left\{ {R^2 2^{ - 2\MyJ\,s}  + [C_1 \tau _\infty   + C(\alpha )\tau _{0,\infty } ]
    \tfrac{{\sqrt {2^{\MyJ} \log \log (n)} }}{n} + [C_1 M^{\,2}  + C(\alpha )M^{\,2} ]
    \tfrac{{2^{\MyJ}  \cdot [\log \log (n)]^2 }}{{n^2 }}} \right\} +  \hfill \\
& & \quad\quad\quad
   + [2C(\alpha )(\tau _{0,\infty }  + \tfrac{1}{3}M\|f_0\|_\infty  ) 
   + C_2 ]\tfrac{{\log \log (n)}}{n},
\end{eqnarray*}
where the last inequality denoted by $(\diamondsuit)$ comes from the fact that, for $n \geqslant 16$, $1 \leqslant \log\log(n)$.

Now, since by hypothesis $2^\MyJ \leqslant n^2/[\log\log(n)]^3$, we have
{\scriptsize{
\begin{eqnarray*}
\frac{{2^{\MyJ}  \cdot [\log \log (n)]^2 }}{{n^2 }} &=& 
    \frac{{\sqrt {2^{\MyJ}  \cdot \log \log (n)} }}{n}
    \frac{{\sqrt {2^{\MyJ}  \cdot \log \log (n)} }}{n}
    \log \log (n) \leqslant  \hfill \\
&\leqslant& 
    \frac{{\sqrt {2^{\MyJ}  \cdot \log \log (n)} }}{n} 
    \cdot \frac{n}{n} \cdot \sqrt {\frac{{\log \log (n)}}
    {{[\log \log (n)]^3 }}} \log \log (n) = \frac{{\sqrt {2^{\MyJ}  \cdot \log \log (n)}}}{n},
\end{eqnarray*}
}}
so that
\begin{eqnarray*}
& & \mathop {\inf }\limits_{\MyJ \in \mathcal{J}_n }
    \left\{ {\|f - \Pi _{\text{S}_\MyJ }(f)|_{\Ldue(G_X )}^2  
    + \widetilde{r}_{n,\MyJ} (\alpha ) + \textsf{V}_{n,\MyJ} (\beta )} \right\} \leqslant \\
&\leqslant& 
    C'\mathop {\inf }\limits_{\MyJ \in \mathcal{J}_n } \left\{ {R^2 2^{ - 2 \MyJ\,s}  
    + 2^{\MyJ/2} \frac{{\sqrt {\log \log (n)} }}{n}} \right\} + 
    C''\frac{{\log \log (n)}}{n},
\end{eqnarray*}
where $C' = 2 \cdot \max \{ 1,C_1 \tau _\infty   + C(\alpha )\tau _{0,\infty } ,M^{\,2} [C_1 + C(\alpha )]\}$ and $C'' = [2C(\alpha )(\tau _{0,\infty }  + \tfrac{1}
{3}M\|f_0\|_\infty  ) + C_2 ]$.

From this point on, the proof continues as in \cite{Fromont:Laurent:2005} and it will reported here just for the sake of completeness. First of all notice that 
\[
R^2 2^{ - 2\MyJ\,s}  \leqslant 2^{\MyJ/2} \tfrac{{\sqrt {\log \log (n)} }}{n} 
    \Leftrightarrow 2^{\MyJ}  > \left[ {\tfrac{{(n\,R)^2 }}{{\log \log (n)}}} 
    \right]^{\tfrac{1}{{1 + 4s}}}.
\]
So define $\MyJ^\star$ by
\[
\MyJ^{\star}  = \left\{ {\log _2 \left[ {\tfrac{{(n\,R)^2 }}{{\log \log (n)}}} \right]^{\tfrac{1}{{1 + 4s}}} } \right\} + 1.
\]
Then we distinguish the following three cases:
\begin{enumerate}
\item In this case we work under the hypothesis that $\MyJ^\star \in \mathcal{J}_n$. This means that 
\[
\MyJ^\star \leqslant \overline{\MyJ}_n = \log_2\big\{ \tfrac{n^2}{[\log\log(n)]^3} \big\},
\]
and that
\[
\mathop {\inf }\limits_{\MyJ \in \mathcal{J}_n } \left\{ {R^2 2^{ - 2 \MyJ \,s}  + 
    2^{\MyJ/2} \tfrac{{\sqrt {\log \log (n)} }}{n}} \right\} 
    \leqslant R^2 2^{ - 2\MyJ^{\star} \,s}  + 2^{\MyJ^{\star} /2} 
    \tfrac{{\sqrt {\log \log (n)} }}{n}.
\]
Now notice that
\begin{itemize}
\item 
    $R^2 2^{ - 2\MyJ^{\star} \,s}  \leqslant R^{\tfrac{2}{{4s + 1}}} 
    \biggl[ {\tfrac{{\sqrt {\log \log (n)} }}{n}} \biggl]^{\tfrac{{4s}}{{4s + 1}}},$
\item 
    $2^{\MyJ^{\star}/2} \tfrac{{\sqrt {\log \log (n)} }}{n} \leqslant 
    \sqrt{2} \left[ {\tfrac{{n\,R}}{{\log \log (n)}}} \right]^{\tfrac{1}{{4s + 1}}} 
    \tfrac{{\sqrt {\log \log (n)} }}{n} \leqslant 
    \sqrt{2} R^{\tfrac{2}{{4s + 1}}} 
    \biggl[ {\tfrac{{\sqrt {\log \log (n)} }}{n}} \biggr]^{\tfrac{{4s}}{{4s + 1}}}.$
\end{itemize}
So we can write
\[
\mathop {\inf }\limits_{\MyJ \in \mathcal{J}_n } 
    \left\{ {R^2 2^{ - 2 \MyJ\,s}  + 2^{\MyJ/2} \tfrac{{\sqrt {\log \log (n)} }}{n}} \right\} 
    \leqslant (1 + \sqrt 2 )R^{\tfrac{2}{{4s + 1}}} 
    \biggl[ {\tfrac{{\sqrt {\log \log (n)} }}{n}} \biggr]^{\tfrac{{4s}}{{4s + 1}}} .
\]
\item In this second case, we assume that $\MyJ^\star > \overline{\MyJ}_n$, hence, by definition of $\MyJ^\star$, for all $\MyJ \in \mathcal{J}_n$, 
\[
2^{\MyJ/2} \tfrac{{\sqrt {\log \log (n)} }}{n} \leqslant R^2 2^{ - 2\MyJ\,s}.
\]
Consequently we get
\[
\mathop {\inf }\limits_{\MyJ \in \mathcal{J}_n } 
    \left\{ {R^2 2^{ - 2\MyJ\,s}  + 2^{\MyJ/2} 
    \tfrac{{\sqrt {\log \log (n)} }}{n}} \right\} \leqslant R^2 2^{ - 2\overline {\MyJ} _n \,s}
    \leqslant 2^{2s + 1} R^2 \left\{ {\tfrac{{[\log \log (n)]^3 }}{{n^2 }}} \right\}^{2s} .
\]
\item In this last case, we assume $\MyJ^\star < 0$. Under this hypothesis, by definition of $\MyJ^\star$,
\[
R^2 2^{-2\MyJ\,s}\leqslant 2^{\MyJ/2} \tfrac{{\sqrt {\log \log (n)} }}{n}, \quad \forall \; \MyJ \in \mathcal{J}_n.
\]
Taking $\MyJ \equiv 0$, we get
\[
\mathop {\inf }\limits_{\MyJ \in \mathcal{J}_n } \left\{ {R^2 2^{ - 2\MyJ\,s}  + 2^{\MyJ/2} 
    \tfrac{{\sqrt {\log \log (n)} }}{n}} \right\} 
    \leqslant \tfrac{{\sqrt {\log \log (n)} }}{n}.
\]
\end{enumerate}
And this complete the proof of Corollary \ref{cor:CH2:unifSepRates}.

\subsection{Proof of Lemma \ref{lem:CH2:HoudreReyn_inAction}}\label{subsec:CH2:proofHoudre}

We can prove Lemma \ref{lem:CH2:HoudreReyn_inAction} by  either using Theorem 3.4 in \cite{Houdre:Reynaud:2002}, or Theorem 3.3 in \cite{Gine:Latala:Zinn:2000}. From these results we know that exists some absolute constant $C > 0$ such that, for all $u > 0$,
\[
\Prob_f^{ \otimes n} \left\{ {|\widetilde{U}_{n,\MyJ} | > \frac{C}{{n(n - 1)}}
    \left[ {\textsf{A}_1 \sqrt u  + \textsf{A}_2 u + \textsf{A}_3 u^{\tfrac{3}{2}}  
    + \textsf{A}_4 u^2 } \right]} \right\} \leqslant 5.6\,\mathrm{e}^{ - u}.
\]
where
\begin{itemize}
\item $\textsf{A}^2_1  = n(n - 1) \, \Exp\big[g_{\MyJ}^2 (\bZ_1 ,\bZ_2 )\big],$
\item $\textsf{A}_2  = \sup \left\{ {\; \left| \; {\Exp\biggl[ {\sum\limits_{i \ne j} 
    {g_{\MyJ} (\bZ_1 ,\bZ_2 )\, a_i (\bZ_1 )\, b_j (\bZ_2 )} } \biggr]} \; \right| : 
    \Exp\left[ {\sum\limits_{i = 1}^n {a_i^2 (\bZ_1 )} } \right] \leqslant 1 
    \text{ and } \Exp\left[ {\sum\limits_{j = 1}^n {b_j^2 (\bZ_2 )} } \right] \leqslant 1} \right\},$
\item $\textsf{A}^2_3  = n\mathop {\sup }\limits_{\bz} \left\{ {\Exp\left[ 
    {g_{\MyJ}^2 (\bz,\bZ_2 )} \right]} \right\},$
\item $\textsf{A}_4  = \mathop {\sup }\limits_{\bz_1 ,\bz_2 } \left| 
    {g_{\MyJ} (\bz_1 ,\bz_2 )} \right|,$
\end{itemize}
and, from Equation (\ref{eq:CH2:2ndOrderUstat}),
\[
g_{\MyJ} (\bz_1,\bz_2 ) = \sum\limits_{k = 1}^{\bar k(\MyJ)} {\big\{ y_1 \phi _{\MyJ,k} (G_X(x_1 )) - \theta_{\MyJ,k} \big\} \cdot \big\{ y_2 \phi _{\MyJ,k} (G_X(x_2 )) - \theta_{\MyJ,k} \big\} },
\]
with $\theta_{\MyJ,k} = \left\langle {f,\phi _{\MyJ,k} (G_X )} \right\rangle_{\Ldue(G_X )}  = \Exp\left[ {Y\phi _{\MyJ,k} (G_X (X))} \right].$

In the following we will bound separately each of these four terms using some specific properties of the \emph{warped wavelet} basis introduced in Section \ref{subsec:CH2:TestPower}. In this section, $\phi(\cdot)$ is the compactly supported scaling function used to generate our basis. If $\textsf{supp}(\phi) \subset [0,L]$ then, for any $k$ and $j$ in $\MyZ$ we put
\[
\textsf{I}_{\MyJ,k}  = \left[ {\tfrac{k}{{2^{\MyJ} }},\tfrac{{k + 1}}{{2^{\MyJ} }}} \right]
\quad \text{and} \quad 
\tilde{\textsf{I}}_{\MyJ,k} = \left[ {\tfrac{k}{{2^{\MyJ}}},\tfrac{{k + L}}{{2^{\MyJ}}}}\right],
\]
so that $\textsf{supp}(\phi_{\MyJ,k}) \subset \tilde{\textsf{I}}_{\MyJ,k}$. Notice that
\[
|k_1  - k_2 | > L\; \Rightarrow \;\tilde{\textsf{I}}_{\MyJ,k_1 }  \cap \tilde{\textsf{I}}_{\MyJ,k_2 }  = \emptyset,
\]
and
\[ 
\textsf{supp}(\phi_{\MyJ,k}\circ G_X) \subset G_X^{-1}\big( \tilde{\textsf{I}}_{\MyJ,k }\big) = \textsf{I}^{G}_{\MyJ,k } \quad \text{with} \quad
\mathbb{1}_{\textsf{I}_{\MyJ,k}^G }(x) = 1 \Leftrightarrow \mathbb{1}_{\tilde{I}_{\MyJ,k}}\big(G_X (x)\big) = 1.
\]

\subsubsection{A bound for $\textsf{A}_1$}\label{subsub:CH2:A1}
Since 
\[
g_{\MyJ}^2 (\bz_1 ,\bz_2 ) = \sum\limits_{k,k'} {\{ y_1 \phi _{\MyJ,k} (G(x_1 )) 
    - \theta _{\MyJ,k } \}\{ y_1 \phi _{\MyJ,k'} (G(x_1 )) 
    - \theta _{\MyJ,k'} \}\{ y_2 \phi _{\MyJ,k } (G(x_2 )) 
    - \theta _{\MyJ,k } \}\{ y_2 \phi _{\MyJ,k'} (G(x_2 )) 
    - \theta _{\MyJ,k'} \} },
\]
by the independence and identical distribution of the sample $\{\bZ_i\}_{i\in\nSeq}$, we have
\[
\Exp\left[ {g_{\MyJ}^2 (\bZ_1 ,\bZ_2 )} \right] = 
    \sum\limits_{k,k'} {\left\{ {\Exp {\big( {Y\phi_{\MyJ,k} (G(X)) 
    - \theta _{\MyJ,k} } \big)\big( {Y\phi_{\MyJ,k'} (G(X)) - \theta _{\MyJ,k'} } \big)}     } \right\}^2 }. 
\]
Now
\be
& & \Exp\left( {Y\phi_{\MyJ,k}(G(X)) - \theta_{\MyJ,k} } \right)
    \left( {Y\phi_{\MyJ,k'}(G(X)) - \theta_{\MyJ,k'} } \right) =  \hfill \\
&=& \Exp\left[ {(f(X) + \epsilon )\phi_{\MyJ,k}(G(X)) 
    - \theta_{\MyJ,k} } \right] \cdot \left[ {(f(X) + \epsilon )\phi _{\MyJ,k'} (G(X)) 
    - \theta _{\MyJ,k'} } \right] =  \hfill \\
&=& \Exp\left[ {f^2 (X)\phi_{\MyJ,k} (G(X))\phi_{\MyJ,k'} (G(X))} \right] 
    - \theta_{\MyJ,k'} \Exp[f(X)\phi_{\MyJ,k} (G(X))] 
    - \theta_{\MyJ,k} \Exp[f(X)\phi _{\MyJ,k'} (G(X))] +  \hfill \\
& & \quad
    + \theta_{\MyJ,k} \theta_{\MyJ,k'} + \Exp\left[ {\phi _{\MyJ,k} (G(X))
    \phi _{\MyJ,k'}(G(X)) \Exp(\epsilon ^2 |X)} \right] =  \hfill \\
&=& \Exp\left\{ {[f^2(X) + \sigma ^2 (X)]\phi _{\MyJ,k} (G(X))\phi _{\MyJ,k'} (G(X))} \right\} - \theta _{\MyJ,k} \theta _{\MyJ,k'}.
\ee
Hence, defining $\tau(x) = f^2(x) + \sigma^2(x)$ and using the inequality $(a-b)^2 \leqslant 2(a^2 + b^2)$ we get 
\be
\Exp\left[ {g_{\MyJ}^2 (\bZ_1 ,\bZ_2 )} \right] &=&
  \sum\limits_{k,k'} {\left\{ {\Exp\left\{ {\tau (X)\phi _{\MyJ,k} (G(X))
  \phi_{\MyJ,k'} (G(X))} \right\} - \theta_{\MyJ,k} \theta _{\MyJ,k'} } \right\}^2 }  
  \leqslant \hfill \\
&\leqslant& 2\sum\limits_{k,k'} {\left\{ {\Exp\left[ {\tau (X)\phi _{\MyJ,k} (G(X))
    \phi _{\MyJ,k'} (G(X))} \right]} \right\}^2 }  + 2\sum\limits_{k,k'} 
    {(\theta _{\MyJ,k} \theta _{\MyJ,k'} )^2 } \leqslant \hfill \\
&\leqslant& 2\sum\limits_{k,k'} {\left\{ {\Exp\left[ {\tau (X)\phi _{\MyJ,k} (G(X))
    \phi _{\MyJ,k'} (G(X))} \right]} \right\}^2 }  + 2 \biggl( {\sum\limits_k 
    {\theta _{\MyJ,k}^2 } } \biggr)^2.
\ee
At this point we proceed bounding separately the two terms in the previous equation.
\begin{itemize}
\item Let
\[
\mathscr{E}_{k,k'}  = \textsf{I}_{\MyJ,k}^G \cap \textsf{I}_{\MyJ,k'}^G  
\Rightarrow \mathscr{E}_{k,k'} \subset \textsf{I}_{\MyJ,k}^G 
\;\text{and}\;
\mathscr{E}_{k,k'}  \subset \textsf{I}_{\MyJ,k'}^G,
\]
and 
\[
\mathscr{I}(k) = \big\{\ell\in\MyZ: \tilde{\textsf{I}}_{\MyJ,k} \cap \tilde{\textsf{I}}_{\MyJ,\ell} \neq \emptyset\big\} = \big\{\ell\in\MyZ:|k-\ell|\leqslant L\big\},
\]
with $\textsf{card}(\mathscr{I}(k)) = 2L + 1$. Hence
\begin{eqnarray}\label{eq:CH2:A1_A}
& & \sum\limits_{k,k'} {\left\{ {\Exp\left[ {\tau (X)\phi _{\MyJ,k} (G(X))
    \phi _{\MyJ,k'} (G(X))} \right]} \right\}^2 } \leqslant\big[2^{\MyJ}\|\phi\|_\infty^2\big]^2 \sum\limits_{k,k'} {{\left\{ {\Exp\left[ 
    {\tau (X)\mathbb{1}_{\mathscr{E}_{k,k'} } (X)} \right]} \right\}^2 } }  = 
    \nonumber \hfill \\
&=& \big[2^{\MyJ}\|\phi\|_\infty ^2 \big]^2 \sum\limits_{k,k'} 
    {{\Exp\left[{\tau (X)\mathbb{1}_{\mathscr{E}_{k,k'}}(X)} 
    \right] \Exp\left[ {\tau (X) \mathbb{1}_{\mathscr{E}_{k,k'}}(X)} \right]} } 
    \leqslant \nonumber \hfill \\
&\leqslant& \big[2^{\MyJ}\|\phi\|_\infty ^2\big]^2 \sum\limits_{k,k'} { 
    {\Exp\left[ {\tau    (X)\mathbb{1}_{\textsf{I}_{\MyJ,k}^G } (X)} \right]
    \Exp\left[ {\tau(X) \mathbb{1}_{\textsf{I}_{\MyJ,k'}^G } (X)} \right]}}
    \leqslant \nonumber \hfill \\
&\leqslant& \big[\tau _\infty  2^{\MyJ}\|\phi\|_\infty ^2\big]^2 \sum\limits_{k,k'} { {\Exp\left[ {\mathbb{1}_{\tilde{\textsf{I}}_{\MyJ,k} } \big(G_X (X)\big)} \right]
    {\Exp\left[ {\mathbb{1}_{\tilde{\textsf{I}}_{\MyJ,k'} }
    \big(G_X (X)\big)} \right]} }} = \nonumber \hfill \\
&=& \big[\tau _\infty  2^{\MyJ}\|\phi\|_\infty ^2\big]^2 \sum\limits_k {\left\{
    {\int_{[0,1]}{\mathbb{1}_{\tilde{\textsf{I}}_{\MyJ,k}}(x)\mathrm{d}x} \sum\limits_{k' \in \mathscr{I}(k)}{\int_{[0,1]}{\mathbb{1}_{\tilde{\textsf{I}}_{\MyJ,k'}}(x)\mathrm{d}x}} } \right\}} \leqslant \nonumber \hfill \\
&\leqslant& \big[\tau _\infty  2^{\MyJ}\|\phi\|_\infty ^2\big]^2 \frac{L(2L+1)}{2^{\MyJ}}\sum\limits_k {\left\{
    {\int_{\tilde{\textsf{I}}_{\MyJ,k}}{\mathbb{1}_{[0,1]}(x)\mathrm{d}x}}\right\}} \leqslant
    2^{\MyJ} \, \big[\tau _\infty \|\phi\|_\infty ^2\big]^2 L^2(2L+1).
\end{eqnarray}
where the last inequality comes from the fact that for any function $w\in\Lp^2(\MyR)$
\[
\sum\limits_{k\in\MyZ} {\left\{ {\int_{\tilde{\textsf{I}}_{\MyJ,k} } {w(x)\mathrm{d}x} } \right\}} \leqslant L\int {w(x)\mathrm{d}x}.
\]
\item We have the following two bounds
\begin{enumerate}
\item $\sum\limits_k {\theta_{\MyJ,k}^2 } \leqslant 
    \sum\limits_{\MyJ}{\sum\limits_k {\theta_{\MyJ,k}^2}} = \|f\|_{\Ldue(G_X )}^2  
    \leqslant \|f\|_\infty ^2.$
\item By using again the inequality we just mentioned, we obtain
\be
\sum\limits_k {\theta _{\MyJ,k}^2 }  &=& \sum\limits_k {\left\{ 
    {\Exp\left[ {Y\phi _{\MyJ,k} (G_X (X))} \right]} \right\}^2 }  = 
    \sum\limits_k {\left\{{\Exp\left[ {f(X)\phi_{\MyJ,k} (G_X (X))} \right]} \right\}^2 } \leqslant  \hfill \\
&\leqslant& 2^{\MyJ} \, \|\phi\|_\infty ^2 \|f\|_\infty  
    \sum\limits_k{\left\{ {\int_{\tilde{\textsf{I}}_{\MyJ,k} } {f(G_X^{ - 1} (x))\mathrm{d}x} }
    \right\}^2 }  \leqslant 2^{\MyJ} \, \|\phi\|_\infty ^2 \|f\|_\infty L\int_{[0,1]} 
    {f(G_X^{ - 1} (x))\mathrm{d}x}  \leqslant  \hfill \\
&\leqslant& 2^{\MyJ} \, \|f\|_\infty ^2 \|\phi\|_\infty ^2 L.
\ee
\end{enumerate}
Combining these two inequalities we can write
\bequ\label{eq:CH2:A1_B}
\biggl( {\sum\limits_k {\theta _{\MyJ,k}^2 } } \biggr)^2  \leqslant 2^{\MyJ} \|f\|_\infty ^4 \|\phi\|_\infty ^2 L \leqslant 2^{\MyJ} \tau^2_\infty \|\phi\|^2_\infty L.
\eequ
\end{itemize}
Finally, from Equations (\ref{eq:CH2:A1_A}) and (\ref{eq:CH2:A1_B}), we obtain
\bequ\label{eq:CH2:A1_final}
\textsf{A}^2_1 \leqslant n(n-1) \, 2^{\MyJ} \tau^2_\infty C'_1(\phi) \Rightarrow 
\textsf{A}_1 \leqslant n \, C_1(\phi) \, \sqrt{2^{\MyJ} \tau^2_\infty}.
\eequ
where $C^2_1(\phi) = 2 \|\phi\|^2_\infty L\{1+\|\phi\|^2_\infty L (2L + 1)\}$.

\subsubsection{A bound for $\textsf{A}_2$}\label{subsub:CH2:A2}
Let $\varphi_{\MyJ,k}(\bz) = y\phi_{\MyJ,k}(G_X(x))$, then
\be
& & \Exp\left\{ {g_{\MyJ} (\bZ_1 ,\bZ_2 )a_i (\bZ_1 )b_j (\bZ_2 )} \right\} = 
    \sum\limits_k {\big\{ {\Exp\left[ {a_i (\bZ_1 )\left( {\varphi _{\MyJ,k} (\bZ_1 ) 
    - \theta _{\MyJ,k} } \right)} \right]} \big\}\big\{ {\Exp[ {b_j (\bZ_2 )
    \left( {\varphi_{\MyJ,k} (\bZ_2 ) - \theta _{\MyJ,k} } \right)}]} 
    \big\}}  =  \hfill \\
&=& \sum\limits_k {\biggl\{ {\Exp[ {a_i (\bZ)\varphi _{\MyJ,k} (\bZ)}] 
    \cdot \Exp[ {b_j (\bZ)\varphi _{\MyJ,k} (\bZ)} ] - \theta _{\MyJ,k} 
    \Exp[ {a_i(\bZ)\varphi_{\MyJ,k}(\bZ)} ] \cdot \Exp[ {b_j (\bZ)}] + }     \biggl.}  \hfill \\
& & \quad\quad \biggl. { - \theta _{\MyJ,k} \Exp [ {b_j (\bZ)\varphi _{\MyJ,k} (\bZ)}] \cdot \Exp [ {a_i (\bZ)} ] + \theta _{\MyJ,k}^2 \Exp[ {a_i (\bZ)}]
    \Exp [ {b_j (\bZ)} ]} \biggl\} =  \hfill \\
&=& \left\{ {\sum\limits_k {\Exp [ {a_i (\bZ)\varphi _{\MyJ,k} (\bZ)} ] 
    \cdot \Exp [ {b_j(\bZ)\varphi_{\MyJ,k}(\bZ)} ]} } \right\} + 
    \left\{ {\Exp\left[ {a_i (\bZ)} \right] \Exp [ {b_j (\bZ)} ]
    \sum\limits_k {\theta_{\MyJ,k}^2 } } \right\} +  \hfill \\
& &\quad\quad - \left\{ {\Exp\left[ {a_i (\bZ)} \right] \cdot \Exp\biggl[ 
    {b_j (\bZ)\sum\limits_k {\theta _{\MyJ,k} \varphi _{\MyJ,k} (\bZ)} } \biggr]} \right\} 
    - \left\{ {\Exp [ {b_j (\bZ)} ] \cdot \Exp\biggl[ {a_i (\bZ)\sum\limits_k 
    {\theta _{\MyJ,k} \varphi _{\MyJ,k} (\bZ)} } \biggr]} \right\} = \hfill \\
&=&   \big\{ \; \textsf{I} \; \big\} + \big\{ \; \textsf{II} \; \big\} 
    - \big\{ \; \textsf{III} \; \big\} - \big\{ \; \textsf{IV} \; \big\}.
\ee
Notice that
\[
\sum\limits_k {\theta _{\MyJ,k} \varphi _{\MyJ,k} (\bZ)}  = 
    \sum\limits_k {\theta _{\MyJ,k} Y\phi _{\MyJ,k} (X)}  = 
    Y\left\{ {\sum\limits_k {\theta _{\MyJ,k} \phi _{\MyJ,k} (X)} } \right\} =
    Y\,\Pi _{\text{S}_\MyJ } (f)(X).
\]
Next, we will bound separately each of the four terms in the previous Equation.
\begin{enumerate}
\item By applying the Cauchy--Schwarz inequality twice we get
\be
\left| {\;\textsf{I}\;} \right| &\leqslant& \left\{ {\sum\limits_k 
    {\Exp\left[ {a_i (\bZ)\varphi _{\MyJ,k}(\bZ) \mathbb{1}_{\mathscr{Z}(k)}(\bZ)} \right]^2 } }
    \right\}^{\tfrac{1}{2}} \left\{ {\sum\limits_k {\Exp\left[ {b_j(\bZ)\varphi_{\MyJ,k}(\bZ)
    \mathbb{1}_{\mathscr{Z}(k)}(\bZ)} \right]^2 } } \right\}^{\tfrac{1}{2}}  \leqslant  \hfill \\
&\leqslant& \left\{ {\sum\limits_k {\Exp\left[ {a_i^2 (\bZ)\mathbb{1}_{\mathscr{Z}(k)}(\bZ)} \right] \Exp\left[ {\varphi _{\MyJ,k}^2 (\bZ)} \right]} } \right\}^{\tfrac{1}{2}} 
    \left\{ {\sum\limits_k {\Exp\left[ {b_j^2(\bZ)\mathbb{1}_{\mathscr{Z}(k)}(\bZ)} 
    \right]\Exp\left[ {\varphi_{\MyJ,k}^2 (\bZ)} \right]} } \right\}^{\tfrac{1}{2}},
\ee
where $\mathscr{Z}(k) = \big\{\bz=(x,y)\in(0,1)\times[-M,M] : x \in \textsf{I}^G_{\MyJ,k}\text{ and } y \in [-M,M] \big\}$. Now we have
\begin{itemize}
\item By usual arguments
\be
\Exp\left[ {\varphi _{\MyJ,k}^2 (\bZ)} \right] &=& \Exp\left[ {Y^2 \phi _{\MyJ,k}^2 (X)} \right] = \Exp\left[ {\tau (X)\phi _{\MyJ,k}^2 (X)} \right] \leqslant \tau _\infty  2^{\MyJ} \,
    \|\phi\|_\infty ^2 \int {\mathbb{1}_{\tilde{\textsf{I}}_{\MyJ,k} } (G_X (x))
    \mathrm{d}G_X (x)}  \leqslant  \hfill \\
&\leqslant& \tau _\infty  2^{\MyJ} \, \|\phi\|_\infty ^2 \tfrac{L}{{2^{\MyJ} }} = 
    \tau _\infty \|\phi\|_\infty ^2 L. \hfill \\ 
\ee
\item $\sum\limits_k {\Exp\left[ {a_i^2 (\bZ)\mathbb{1}_{\mathscr{Z}(k)}(\bZ)} \right]} \leqslant \Exp\left[ {a_i^2 (\bZ)} \right]\sum\limits_k {\Exp\left[ 
    {\mathbb{1}_{\mathscr{Z}(k)}(\bZ)} \right]} \leqslant L\,
    \Exp\left[ {a_i^2 (\bZ)} \right].$
\item $\sum\limits_k {\Exp\left[ {b_j^2 (\bZ)\mathbb{1}_{\mathscr{Z}(k)}(\bZ)} \right]} \leqslant \Exp\left[ {b_j^2 (\bZ)} \right]\sum\limits_k {\Exp\left[ 
    {\mathbb{1}_{\mathscr{Z}(k)}(\bZ)} \right]} \leqslant L\,
    \Exp\left[ {b_j^2 (\bZ)} \right].$
\end{itemize}
Consequently we get
\[
\left| {\;\textsf{I}\;} \right| \leqslant \left\{ {\tau _\infty  \|\phi\|_\infty ^2 L^2 \Exp\left[ {a_i^2 (\bZ)} \right]} \right\}^{\tfrac{1}{2}} 
    \left\{ {\tau _\infty \|\phi\|_\infty ^2 L^2 
    \Exp\left[ {b_j^2 (\bZ)} \right]} \right\}^{\tfrac{1}{2}}  
    = \tau _\infty \|\phi\|_\infty ^2 L^2 \sqrt {\Exp[ {a_i^2 (\bZ)}]} 
    \sqrt {\Exp[ {b_j^2 (\bZ)}]}. 
\]
And finally
\be
\sum\limits_{i = 1}^n {\sum\limits_{j = 1}^{i - 1} {\left| {\;\textsf{I}\;} \right|} } &\leqslant& \tau _\infty  |\phi |_\infty ^2 L^2 \left\{ 
    {\sum\limits_{i = 1}^n {\sqrt {\Exp[a_i^2 (\bZ)]} } \; \cdot \;
    \sum\limits_{j = 1}^{i - 1} {\sqrt {\Exp[b_j^2 (\bZ)]} } } \right\} \leqslant \hfill \\
&\leqslant& n\;\tau _\infty  \|\phi\|_\infty ^2 L^2 
    \left\{ {\sqrt {\sum\limits_{i = 1}^n {\Exp[a_i^2 (\bZ)]} } \; \cdot \;
    \sqrt {\sum\limits_{j = 1}^n {\Exp[b_j^2(\bZ)]} } } \right\} 
    \leqslant n\;\tau _\infty  \|\phi \|_\infty ^2 L^2.
\ee
\item By definition of $a_i(\cdot)$ and $b_j(\cdot)$, we obtain
\be
\sum\limits_{i = 1}^n {\sum\limits_{j = 1}^{i - 1} {\left| {\;\textsf{II}\;} \right|} } &\leqslant& \sum\limits_{i = 1}^n {\sum\limits_{j = 1}^{i - 1} {\biggl( 
    {\sum\limits_k {\theta _{\MyJ,k}^2 } } \biggr)
    \left| {\Exp[a_i^2 (\bZ)] \cdot \Exp[b_j^2 (\bZ)]} \right|} }  \leqslant \|f\|_\infty ^2 \sum\limits_{i = 1}^n {\Exp[a_i^2 (\bZ)]} \;\sum\limits_{i = 1}^n 
    {\Exp[a_i^2 (\bZ)]} \leqslant \hfill \\
&\leqslant& n\,\|f\|_\infty ^2 \sqrt {\Exp\sum\nolimits_i {a_i^2 (\bZ)} } \;\sqrt {\Exp\sum\nolimits_j {b_j^2 (\bZ)} } \leqslant n\,\|f\|_\infty ^2.
\ee
\item By Cauchy--Schwarz inequality we have
\be
\sum\limits_{i = 1}^n {\sum\limits_{j = 1}^{i - 1} {\left| {\;\textsf{III}\;} \right|} } &\leqslant& \sum\limits_{i = 1}^n {\sum\limits_{j = 1}^{i - 1} 
    {\left| {\Exp[a_i (\bZ)] \cdot \Exp[b_j (\bZ)\,Y\,\Pi_{\text{S}_{\MyJ} }(f)(X)]} \right|} }
    \leqslant \hfill \\
&\leqslant& \sum\limits_{i = 1}^n {\sum\limits_{j = 1}^{i - 1} {\left\{ {
    \Exp[a_i^2 (\bZ)]} \right\}^{\tfrac{1}{2}} \left\{ {\Exp[b_j^2 (\bZ)]} 
    \right\}^{\tfrac{1}{2}} \left\{ {\Exp[Y\,\Pi _{\text{S}_{\MyJ} } (f)(X)]^2 } 
    \right\}^{\tfrac{1}{2}} } }  \leqslant  \hfill \\
&\leqslant& n\left\{ {\Exp[Y{}^2\,\Pi_{\text{S}_{\MyJ} }^2 (f)(X)]} 
    \right\}^{\tfrac{1}{2}} \sqrt {\Exp\sum\nolimits_i {a_i^2 (\bZ)} } \;
    \sqrt {\Exp\sum\nolimits_j {b_j^2(\bZ)} } \leqslant n 
    \left\{ {\Exp[Y{}^2\,\Pi _{\text{S}_{\MyJ} }^2 (f)(X)]} \right\}^{\tfrac{1}{2}}  
    \leqslant  \hfill \\
&\leqslant& n\,\Exp\left\{ {\tau (X)\,\Pi _{\text{S}_{\MyJ} }^2 (f)(X)} 
    \right\}^{\tfrac{1}{2}} \leqslant n\,\sqrt{\tau _\infty \Exp [ {f^2 (X)} ]}
    \leqslant n\,\sqrt{\tau _\infty  \|f\|_\infty ^2} \leqslant n\,\tau _\infty.
\ee
\item Proceeding as in the previous point, we get
\[
\sum\limits_{i = 1}^n {\sum\limits_{j = 1}^{i - 1} {\left| {\;\textsf{IV}\;} \right|} } \leqslant \sum\limits_{i = 1}^n {\sum\limits_{j = 1}^{i - 1} {\left| {\Exp[b_j (\bZ)] \cdot \Exp[a_i(\bZ)\,Y\,\Pi _{\text{S}_{\MyJ} } (f)(X)]} \right|} }  \leqslant n\,\tau _\infty .
\]
\end{enumerate}
Combining all the previous inequalities, we can write
\bequ\label{eq:CH2:A2_final}
\textsf{A}_2 \leqslant \left\{ {n\,\tau _\infty  \|\phi \|_\infty ^2 L^2  + 3n\,\tau _\infty  } \right\} = C_2 (\phi )\,n\,\tau _\infty  ,
\eequ
where $C_2 (\phi ) = \|\phi\|_\infty ^2 L^2  + 3.$

\subsubsection{A bound for $\textsf{A}_3$}\label{subsub:CH2:A3}

Lets start writing
\be
\Exp\left\{ {g_{\MyJ}^2 (\bz,\bZ_2 )} \right\} &=& 
    \sum\limits_{k,k'} {\big\{ y\phi _{\MyJ,k} (G(x)) - \theta _{\MyJ,k} \big\} 
    \big\{ y\phi _{\MyJ,k'} (G(x)) - \theta _{\MyJ,k'} \big\}  \times }  \hfill \\
& & \quad\quad 
   \times \Exp\big\{ Y_2 \phi _{\MyJ,k} (G(X_2 )) - \theta _{\MyJ,k} \big\} \,\Exp\big\{ Y_2 
   \phi _{\MyJ,k'} (G(X_2)) - \theta _{\MyJ,k'} \big\}  =  \hfill \\
&=& \sum\limits_{k,k'} {\{ \; \cdot \;\} \{ \; \cdot \;\} \left\{ {\Exp\left[ {\tau (X)
    \phi _{\MyJ,k} (G(X))\phi _{\MyJ,k'} (G(X))} \right] - \theta _{\MyJ,k} 
    \theta _{\MyJ,k'} } \right\}}  = \hfill \\
&=& \sum\limits_{k,k'} {\{ \; \cdot \;\} \{ \; \cdot \;\} \,\Exp\left[ {\; \cdot \;} \right]} - \biggl[ {\sum\limits_k {\theta _{\MyJ,k} \left( {y\phi _{\MyJ,k} (G(x)) - \theta _{\MyJ,k} } \right)} } \biggl]^2  \leqslant  \hfill \\
&\leqslant& \sum\limits_{k,k'} {\big\{ y\phi _{\MyJ,k} (G(x)) - \theta _{\MyJ,k} \big\} 
    \big\{ y \phi _{\MyJ,k'} (G(x)) - \theta _{\MyJ,k'} \big\} \,\Exp\left[ {\tau (X)
    \phi _{\MyJ,k} (G(X)) \phi _{\MyJ,k'} (G(X))} \right]} .
\ee
Now we have
\begin{itemize}
\item By the same arguments used in bounding $\textsf{A}_1$, we get
\be
\left| {\,\Exp\left[ {\tau (X)\phi _{\MyJ,k} (G(X))\phi _{\MyJ,k'} (G(X))} \right]\,} \right| &\leqslant& \tau _\infty \big[2^{\MyJ} \, \|\phi \|_\infty ^2 \big]
    \int {\mathbb{1}_{\tilde{\textsf{I}}_{\MyJ,k}  \cap \tilde{\textsf{I}}_{\MyJ,k'} } 
    (G_X (x))\mathrm{d}G_X (x)} = \hfill \\
&=& \tau _\infty  \big[2^{\MyJ} \, \|\phi\|_\infty ^2 \big]
    \int {\mathbb{1}_{\tilde{\textsf{I}}_{\MyJ,k} 
    \cap \tilde{\textsf{I}}_{\MyJ,k'} }(x)\mathrm{d}x} .
\ee
\item Since
\[
\sum\limits_{k,k'} {\big\{ y\phi _{\MyJ,k} (G(x)) - \theta _{\MyJ,k} \big\} 
    \big\{ y\phi _{\MyJ,k'} (G(x)) - \theta _{\MyJ,k'} \} }  = 
    \biggl[ {\sum\limits_k {\big\{ y\phi _{\MyJ,k} (G(x)) - \theta _{\MyJ,k} \big\} } } \biggr]^2 
\]
we need to bound separately $\left|{\sum\nolimits_k{\theta_{\MyJ,k}}}\right|$, and
$\mathop {\sup }\limits_{x,y} \big|\sum\nolimits_k {y\phi _{\MyJ,k} (G(x))}\big|$ as follow:
{\footnotesize{
\begin{enumerate}
\item $\mathop {\sup }\limits_{x,y} \big|\sum\nolimits_k {y\phi _{\MyJ,k} (G(x))} \big| \leqslant M \, 2^{\MyJ/2} \|\phi\|_\infty  (2L + 1),$
\item $\left| {\sum\nolimits_k {\theta_{\MyJ,k} } } \right| = \left| 
    {\Exp\left\{ {f(X)\sum\nolimits_k {\phi _{\MyJ,k'} (G(X))} } \right\}} \right| 
    \leqslant \Exp[f(X)]\,\|\sum\nolimits_k {\phi _{\MyJ,k'}(G)} \|_\infty   
    \leqslant \|f\|_\infty 2^{\MyJ/2} \|\phi \|_\infty  (2L + 1)$.
\end{enumerate}
}}
Consequently
\be
& & \mathop {\sup }\limits_{x,y} \left| {\,\biggl[ {\sum\limits_k {
    \big\{ y\phi _{\MyJ,k} (G(x)) - \theta _{\MyJ,k} \big\} } } \biggr]^2 } \right| 
    \leqslant \left[ {\mathop {\sup }\limits_{x,y} \biggl|\sum\nolimits_k {
    y\phi_{\MyJ,k} (G(x))} - \theta _{\MyJ,k} \biggr|  } \right]^2  \leqslant \hfill \\
&\leqslant& \left[ {\mathop {\sup }\limits_{x,y} 
    \biggl|\sum\nolimits_k {y\phi _{\MyJ,k} (G(x))} \biggr| + \left| 
    {\sum\nolimits_k {\theta _{\MyJ,k} } } \right|} \right]^2  \leqslant \left[ {(M + \|f\|_\infty )2^{\MyJ/2}\|\phi\|_\infty  (2L + 1)} \right]^2.
\ee
\end{itemize}
Hence, finally
\be
\mathop {\sup }\limits_{\bz} \Exp\left\{ {g_{\MyJ}^2 (\bz,\bZ_2 )} \right\} 
&\leqslant& \tau _\infty 2^{\MyJ} \|\phi \|_\infty ^2 \frac{L}{{2^{\MyJ} }}\;
    \mathop {\sup }\limits_{x,y} \left| {\,\biggl[ {\sum\limits_k {\{ y\phi _{\MyJ,k} (G(x)) - \theta _{\MyJ,k} \} } } \biggr]^2 } \right| \leqslant  \hfill \\
&\leqslant& \tau _\infty  (M + \|f\|_\infty  )^2 2^{\MyJ} \|\phi\|_\infty ^2 L
    \left[ {\|\phi\|_\infty (2L + 1)} \right]^2,
\ee
so that
\bequ\label{eq:CH2:A3_final}
\textsf{A}_3 \leqslant C_3(\phi) \,(M + \|f\|_\infty) \, \sqrt{n \tau_\infty 2^{\MyJ}}
\leqslant 2 C_3(\phi) \, M  \, \sqrt{n \tau_\infty 2^{\MyJ}},
\eequ 
with $C_3(\phi) = \|\phi\|^2_\infty \sqrt{L}(2L + 1)$.

\subsubsection{A bound for $\textsf{A}_4$}\label{subsub:CH2:A4}

Bearing in mind the following inequalities
\begin{itemize}
\item $\mathop {\sup }\limits_{x,y} \big|\sum\nolimits_k {y\phi _{\MyJ,k} (G(x))}  
    - \theta _{\MyJ,k} \big| \leqslant \|\phi\|_\infty (2L + 1)(M + \|f\|_\infty  )2^{\MyJ/2}$,
\item $\mathop {\sup }\limits_{x,y} \big|y\phi _{\MyJ,k} (G(x))\big| 
    \leqslant M\,2^{\MyJ/2} \|\phi \|_\infty$,
\item $ |\theta _{\MyJ,k} | \leqslant \Exp\big|f(X)\phi _{\MyJ,k} (G(X))\big| 
    \leqslant \|f\|_\infty 2^{\MyJ/2} \|\phi \|_\infty$,
\end{itemize}
we end up with the following bound
\[
\textsf{A}_4  = \mathop {\sup }\limits_{\bz_1 ,\bz_2 } \left| {g_{\MyJ} (\bz_1 ,\bz_2 )} \right| \leqslant C_4 (\phi )2^{\MyJ} (M + \|f\|_\infty  )^2,
\]
where $C_4(\phi) = \|\phi\|^2_\infty (2L+1)$.

\subsubsection{Conclusion}\label{sec:CH2:Houdre_inAction_Conclude}
Up to now we have found that, for each $u > 0$, $\Prob_f^{ \otimes n} \left\{ {|\widetilde{U}_{n,\MyJ} | > \eta_\MyJ (u)} \right\} \leqslant 5.6 \mathrm{e}^{-u}$, with
\[
\eta_\MyJ (u) = \frac{C}{{n - 1}}\left\{ 
    {\frac{{C_1 (\phi )\,n\,\tau _\infty  \sqrt {2^{\MyJ} } }}{n}u^{\tfrac{1}{2}}  
  + \frac{{C_2 (\phi )\,n\,\tau _\infty  }}{n}u + 
  \frac{2 {C_3 (\phi )M\sqrt {n\,\tau _\infty  \,2^{\MyJ} } }}{n}u^{\tfrac{1}{2}} u 
  + \frac{{C_4 (\phi )M^{\, 2} 2^{\MyJ} }}{n}u^2 } \right\}.
\]
By applying the inequality $2 a b \leqslant a^2 + b^2$, we have
\[
C_3 (\phi )\left\{ {2\biggr[ {\sqrt {\tau _\infty  u} } \biggl] \cdot \biggr[ 
    {\sqrt {\tfrac{{M^{\,2} 2^{\MyJ} }}{n}} u} \biggl]} \right\} 
    \leqslant C_3 (\phi )\left\{ {\tau _\infty u + M^{\,2} 2^{\MyJ} \tfrac{{u^2 }}{n}} \right\},
\]
so
\be
\eta_\MyJ (u) &=& \frac{C}{{n - 1}}\left\{ {C_1 (\phi )\,\tau _\infty  \sqrt {2^{\MyJ} } u^{\tfrac{1}{2}}  + \left[ {C_2 (\phi ) + C_3 (\phi )} \right]\,\tau _\infty  
    u + \left[ {C_3 (\phi ) + C_4 (\phi )} \right]M^2 2^{MyJ} \tfrac{{u^2 }}{n}} \right\} 
    \leqslant  \hfill \\
&\leqslant& \frac{\kappa_0} {{n - 1}}\left\{ {\tau _\infty  \sqrt {2^{\MyJ} u}
    + \tau _\infty  u + M^2 \frac{{2^{\MyJ} u^2 }}{n}} \right\},
\ee
where $\kappa_0 = C\max \left\{ {C_1 (\phi ),C_2 (\phi ) + C_3 (\phi ),C_3 (\phi ) + C_4 (\phi )} \right\}$. And this complete the proof.
\EndTh
\bibliographystyle{plain}
\bibliography{papero_A.bbl}
\end{document}